\documentclass[11pt]{article}
\usepackage[utf8]{inputenc} 
\usepackage[T1]{fontenc}    
\usepackage{amsmath,amssymb,amsfonts}
\usepackage{amsthm}
\usepackage{graphicx}
\usepackage{algorithm,algorithmic}
\usepackage{textcomp}
\usepackage{mathtools}
\usepackage[margin=1in]{geometry}
\usepackage{xcolor}
\usepackage{microtype}      
\usepackage{tikz}
\usetikzlibrary{calc}
\usetikzlibrary{positioning}
\usetikzlibrary{arrows.meta}

\usepackage[font={small}]{caption}
\usepackage{setspace}
\usepackage{url}
\usepackage{here}
\usepackage{booktabs}
\usepackage{multirow}
\usepackage{multicol}
\usepackage{bigstrut}
\usepackage{cases}
\usepackage{comment}
\usepackage{cite}
\usepackage[noblocks]{authblk}

\usepackage{framed}
\newtheorem{theorem}{Theorem}

\newtheorem{lemma}{Lemma}
\newtheorem{exam}{Example}
\newtheorem{prop}{Proposition}
\newtheorem{rem}{Remark}
\newtheorem{defi}{Definition}
\newtheorem{assu}{Assumption}

\usepackage{mathtools} 
\DeclareMathOperator*{\argmin}{arg\,min}

\newcommand{\pd}{\partial}

\newcommand{\CC}{\mathcal{C}}
\newcommand{\CN}{\mathcal{N}}
\newcommand{\CD}{\mathcal{D}}
\newcommand{\CE}{\mathcal{E}}

\newcommand{\CQ}{\mathcal{Q}}

\newcommand{\QG}{ { \mathcal{Q}_\mathcal{G}^{\mathrm{all}} } }
\newcommand{\QGactive}{\mathcal{Q}_\mathcal{G}}
\newcommand{\QiGactive}{\mathcal{Q}_\mathcal{G}^{i}}
\newcommand{\QjGactive}{\mathcal{Q}_\mathcal{G}^{j}}
\newcommand{\maxQG}{\mathcal{Q}_\mathcal{G}^\mathrm{max}}

\newcommand{\QiG}{\mathcal{Q}^{i}_\mathcal{G}}

\newcommand{\CG}{\mathcal{G}}

\newcommand{\hl}{\hat{l}}
\newcommand{\Bx}{\mathbf{x}}

\newcommand{\Bv}{\mathbf{v}}

\newcommand{\By}{\mathbf{y}}

\newcommand{\Bw}{\mathbf{w}}

\newcommand{\Bz}{\mathbf{z}}

\newcommand{\Bs}{\mathbf{s}}

\newcommand{\BD}{\mathbf{D}}
\newcommand{\BL}{\mathbf{L}}

\newcommand{\BQ}{\mathbf{Q}}
\newcommand{\BW}{\mathbf{W}}
\newcommand{\BR}{\mathbb{R}}

\newcommand{\diag}{\text{diag}}
\newcommand{\Bzeta}{\boldsymbol\zeta}
\newcommand{\Bxi}{\boldsymbol\xi}

\newcommand{\BPhi}{{\boldsymbol\Phi}}
\newcommand{\CCl}{{\mathcal{C}_l}}
\newcommand{\prox}{\mathrm{prox}}

\newcommand{\qstar}{|\mathcal{Q}_G|}
\newcommand{\bdiag}{\text{blk-diag}}

\title{Distributed Optimization of Clique-Wise Coupled Problems via Three-Operator Splitting
\thanks{
Yuto Watanabe is with the Department of Electrical and Computer Engineering, University of California San Diego, San Diego, CA 92093 USA (email: \texttt{y1watanabe@ucsd.edu}).
Kazunori Sakurama is with the Department of System Innovation, Graduate School of Engineering Science, Osaka University, 1-3, Machikaneyama, Toyonaka, Osaka 560-8531, Japan (email: \texttt{sakurama.kazunori.es@osaka-u.ac.jp}).
This work was partially supported by
the joint project of Kyoto University and Toyota Motor Corporation, titled ``Advanced Mathematical Science for Mobility Society''.}
}
\author{Yuto Watanabe and Kazunori Sakurama}
\date{\today}

\begin{document}
\maketitle
\begin{abstract}
This study explores distributed optimization problems with clique-wise coupling via operator splitting 
and how we can utilize this framework for performance analysis and enhancement.
This framework extends beyond conventional
pairwise coupled problems (e.g., consensus optimization) and is applicable to broader examples.
To this end,
we first introduce a new distributed optimization algorithm by leveraging a clique-based matrix and
the
Davis-Yin splitting (DYS), a versatile three-operator splitting method.
{We then demonstrate that this approach sheds new light on conventional algorithms in the following way:
(i) Existing algorithms (NIDS, Exact diffusion, diffusion, and our previous work) can be derived from our proposed method;
(ii) We present a new mixing matrix based on clique-wise coupling, which surfaces when deriving the NIDS.
We prove its preferable distribution of eigenvalues, enabling fast consensus;
(iii) These observations yield a new linear convergence rate for the NIDS with non-smooth objective functions.
Remarkably our linear rate is first established for the general DYS with a projection for a subspace. This case is not covered by any prior results, to our knowledge.}
Finally, numerical examples showcase the efficacy of our proposed approach.
\end{abstract}

\section{Introduction}\label{sec:introduction}

The last two decades have witnessed the significant advancement of distributed optimization.
In the literature, a huge body of existing studies has been dedicated to \textit{pairwise coupled optimization problems}, where every coupling of variables comprises two agents' decision variables corresponding to the communication path (edge) between the two.
Representative examples are consensus optimization \cite{nedic2009distributed,shi2015proximal,yuan2018exact,li2019decentralized,alghunaim2020decentralized,xu2021distributed}
and formation control \cite{sakurama2022generalized}.
Moreover, so are the problems with globally coupled linear constraints \cite{liang2019distributed} because their dual problems result in pairwise coupled consensus optimization.

\begin{figure}[t]
    \centering
    \includegraphics[width=0.7\columnwidth]{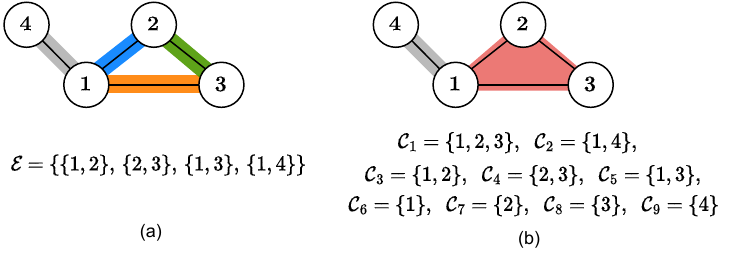}
    \caption{Sketches of (a) pairwise coupling and (b) clique-wise coupling.
    }
    \label{fig:clique_ex}
\end{figure}

To handle wider applications that involve complex coupling beyond edges, we 
leverage \textit{cliques}, complete subgraphs of a graph \cite{bollobas1998modern}, as a generalization of edges and
tackle a more generic class of distributed optimization---\textit{clique-wise coupled optimization problems}.
This class has been introduced in our recent works \cite{watanabe2022distributed,watanabe2022clique-wise} with an emphasis on its generalization aspect.
In this note, we elucidate additional benefits of this class of problems for performance enhancement and analysis via a new algorithm based on a three-operator splitting \cite{davis2017three}.
This class of problems is formulated as follows:
Consider a multi-agent system with $n$ agents over a time-invariant undirected graph $\CG=(\CN,\CE)$ with $\CN=\{1,\ldots,n\}$ and an edge set $\CE$.
Let $x_i\in\BR^{d_i}$ represent the $d_i$ dimensional decision variable of agent $i$.
Then, the following is called a \textit{clique-wise coupled optimization problem}:
\begin{align}\label{problem}
\!\!\!\!\!\!\!\!
\begin{array}{cl}
\underset{
\substack{x_i\in \BR^{d_i}\\\,i\in\CN}
}
{\mbox{min}}&
\!\!\!\!\!
\displaystyle  
 \underbrace{\sum_{l\in\QGactive} \left(f_l (x_\CCl) +g_l(x_\CCl)\right)}_{\text{clique-wise coupling}}
+ \displaystyle  
\sum_{i=1}^n \left(\hat{f}_i (x_i) +\hat{g}_i (x_i) \right),
\end{array}
\end{align}
where the set $\CCl\subset \CN$ represents a clique, and
the set $\QGactive \neq \emptyset$ is a subset of $\QG$, the index set of all the cliques in $\CG$.
(For example, in the undirected graph in Fig. \ref{fig:clique_ex}, we have $\QG=\{1,\ldots,9\}$ and $\CC_1,\ldots,\CC_9$.)
For the set $\CCl=\{j_1,\ldots,j_{|\CCl|}\}\subset\CN$, let $x_{\CCl}$ denote $x_{\CCl}=[x_{j_1}^\top,\ldots,x_{j_{|\CCl|}}^\top]^\top$.
For all $l\in\QGactive$,
$f_l:\BR^{\sum_{j\in\CCl}d_j}\to \BR$ is $L_l$-smooth and convex, and $g_l:\BR^{\sum_{j\in\CCl}d_j}\to \BR$
is proper, closed, and convex, where the subscript "$j$" shows the index of agent $j$ in $\CCl.$
For all $i\in\CN$,
$\hat{f}_i:\BR^{d_i}\to\BR$ is $\hat{L}_i$-smooth and convex, and $\hat{g}_i:\BR^{{d_i}}\to\BR$
is proper, closed, and convex.

As mentioned above,
an immediate benefit of Problem \eqref{problem} is that it can handle variable couplings of more than two agents.
As Fig. \ref{fig:clique_ex},
cliques in (b) can deal with the coupling of three nodes $\{1,2,3\}$, differently from (a).
Indeed, Problem \eqref{problem} always contains conventional pairwise coupled optimization problems as nodes and edges are also cliques.
As well as pairwise coupled problems, 
other possible applications are, for example, (i) clique-wise coupled linear constraints \cite{watanabe2022distributed,watanabe2022clique-wise,li2021primal} (e.g., resource allocation in Section \ref{sec:numerical_experiments}), (ii) sparse SDP \cite{vandenberghe2015chordal} (e.g., distributed design of distributed controllers \cite{zheng2019distributed}, sensor network localization \cite{anjos2011handbook}, etc), (iii) regularization accounting for network structures (e.g, Network lasso \cite{hallac2015network}), and
 (iv) approximation of trace norm minimization problems (e.g., multi-task learning
 \cite{zhang2018overview}, robust PCA \cite{candes2011robust}, etc).

This note addresses Problem \eqref{problem} using the \textit{Davis-Yin Splitting} (DYS) and reveals that the notion of clique-wise coupling is beneficial for analyzing and improving convergence performance. 
The DYS is a versatile three-operator splitting scheme that
generalizes basic operator-splitting methods
(e.g., the forward-backward and Douglas-Rachford splittings).
Firstly,
we reformulate Problem \eqref{problem} by introducing a matrix called the \textit{clique-wise duplication (CD) matrix}.
This matrix lifts Problem \eqref{problem} to a tractable separated form for algorithm design.
Then, applying the DYS, we derive the proposed algorithm called the clique-based distributed DYS (CD-DYS).
Subsequently, we demonstrate that the CD-DYS generalizes several existing algorithms, encompassing the celebrated NIDS \cite{li2019decentralized}.
Then, we analyze a new mixing matrix that naturally comes up in deriving the NIDS and show a preferable distribution of its eigenvalues.
Moreover, we present a new linear convergence rate for the NIDS with non-smooth terms by proving a more general linear rate for the DYS with a projection onto a subspace under strong convexity of the smooth term.
Finally, numerical examples illustrate the effectiveness of the proposed approach.

Our contributions can be summarized as follows.
(i) We propose a new distributed algorithm, CD-DYS, for Problem \eqref{problem} applicable to broad examples ranging from optimization to control and learning problems;
(ii) Our investigation of consensus optimization as a clique-wise coupled problem unveils that several conventional distributed optimization methods, including NIDS \cite{li2019decentralized}, are derived from the proposed CD-DYS method, which leads to a new linear convergence rate for the NIDS with non-smooth objective functions. This linear rate admits bigger stepsizes than ones in \cite{alghunaim2020decentralized,xu2021distributed}.
 It is worth mentioning that our linear convergence is first established for the general DYS with an indicator function of a linear image space, which does not follow from the prior works \cite{davis2017three,lee2022convergence,yi2022convergence,condat2022randprox} as indicator functions are neither smooth nor strongly convex;
(iii) Numerical examples demonstrate the higher performance of our proposed approach than \cite{li2019decentralized} and \cite{liang2019distributed}.
In particular, the superiority against the standard NIDS \cite{li2019decentralized} is attributed to a novel mixing matrix obtained from our proposed method, which realizes a preferable eigenvalue distribution for fast consensus.
We also provide its theoretical evidence.
Note that one can construct this matrix without global information and use it for other consensus-based algorithms.

The remainder of this note is organized as follows.
Section \ref{sec:preliminaries} provides preliminaries.
Section \ref{sec:CDmatrix} presents the definition of the CD matrix and its analysis including a new mixing matrix.
In Section \ref{sec:alg}, we propose new distributed algorithms based on the DYS.
In Section \ref{sec:analysis}, we analyze the proposed methods for consensus optimization and show a new linear convergence result.
Section \ref{sec:numerical_experiments} presents numerical experiments.
Section \ref{sec:proof-rate} provides the proof of the convergence rate.
Finally, Section \ref{sec:conclusion} concludes this note.

\section{Preliminaries}\label{sec:preliminaries}

We here prepare several important notions.
\paragraph*{Notations}
Throughout this note, we use the following notations.
Let $|\cdot|$ be the number of elements in a countable finite set.
Let $I_d$ denote the $d\times d$ identity matrix in $\mathbb{R}^{d\times d}$.
We omit the subscript $d$ of $I_d$ when the dimension is obvious.
Let $O_{d_1\times d_2}$ be the $d_1\times d_2$ zero matrix.
Let $\mathbf{1}_d = [1,\ldots,1]^\top\in \mathbb{R}^{d}$.
For $\mathcal{M}\subset \mathcal{N}$, $[x_j]_{j\in\mathcal{M}}$ and $x_{\mathcal{M}}$ represent the stacked vector in ascending order obtained from vectors $x_j\in\BR^{d_j},\,j\in\mathcal{M}$, and
we use the same notation to express stacked matrices.
Let $\diag(a)$ with $a=[a_1,\ldots,a_n]^\top$ denote the diagonal matrix whose $i$th diagonal entry is $a_i\in\BR$.
Similarly, $\bdiag([\ldots,R_i,\ldots])$ and $\bdiag([R_j]_{j\in\mathcal{M}})$ represent the block diagonal matrix. For a symmetric matrix $Q\succ O$, let $\| u \|_Q=\sqrt{\langle u,u \rangle_Q }$ with the inner product $\langle u,v \rangle_Q := v^\top Q u$, and we simply write $\|\cdot\|_{I_m}=\|\cdot\|$ for $Q=I_m$.
Let $\lambda_\mathrm{max}(Q)$ and $\lambda_\mathrm{min}(Q)$ be the largest and smallest eigenvalues of $Q$, respectively.
For a differentiable function $f:\BR^d\to\mathbb{R}$ and $x\in\BR^d$, we write $\nabla_x f(\cdot)=\pd f/\pd x(\cdot)$.
We simply use $\nabla$ when it is obvious.
For a proper, closed, and convex function $g:\BR^d:\to\BR\cup\{+ \infty\}$ and $Q\succ 0$, 
the proximal operator of $g$ with $Q$ is represented by $\prox^Q_g(x) = \argmin_{x'\in\BR^d} \{ g(x') + \|x-x'\|_Q^2/2\}$, and we denote $\prox_{g}^{I}(\cdot)=\prox_g(\cdot)$ for $Q=I$.
Let $\delta_\CD(\cdot)$ represent the indicator function of $\CD$, i.e., $\delta_\CD(x)=0$ for $x\in\CD$ and $\delta_\CD(x)=\infty$ for $x\notin\CD$.
The projection onto a closed convex set $\CD$ with a metric $Q$ is represented by
$P^Q_\CD(x) = \argmin_{x'\in\CD}\|x-x'\|_Q$, and
we write $P^I_\CD(\cdot)=P_\CD(\cdot)$ for $Q=I$.

\paragraph{Graph theory}\label{subsec:graph_theory}
Consider a graph $\CG=(\CN,\CE)$ with a node set $\CN=\{1,\ldots,n\}$ and an edge set $\CE$ consisting of pairs $\{i,j\}$ of different nodes $i,j\in \CN$.
Note that throughout this note, we consider undirected graphs and do not distinguish $\{i,j\}$ and $\{j,i\}$ for each $\{i,j\}\in\CE$.
For $i\in\CN$ and $\CG$, let $\CN_i\subset \CN$ be the \textit{neighbor set} of node $i$ over $\CG$, defined as $\CN_i=\{j\in \mathcal{N}:\{i,j\}\in \mathcal{E}\}\cup\{i\}$.

For an undirected graph $\CG$,
consider a set $\CC\subset\CN$.
The set $\CC$ is called a \textit{clique} of $\CG$ if the subgraph $\CG$ induced by $\CC$ is complete \cite{bollobas1998modern}.
We define $\QG=\{1,2,\ldots,q\}$ as the set of indices of all the cliques in $\CG$.
For $\QG$, the set $\QGactive$ represents a subset of $\QG$.
If a clique $\CC$ is not contained by any other cliques, $\CC$ is said to be \textit{maximal}.
Let $\maxQG(\subset\QG)$ be the set of indices of all the maximal cliques in $\CG$. 
For edge set $\CE$, let $\QGactive^{\mathrm{edge}}$ be the index set of all the edges.
For $\QGactive\subset\QG$ and $i\in\CN$, we define $\QiGactive$ as the index set of all cliques in $\QGactive$ containing $i$.
Similarly,
$\QGactive^{ij}$ represents 
$\QGactive^{ij}=\QGactive^{ji}=\QGactive^i\cap\QGactive^j$.
For each $i\in\CN$, $\CN_i$, and $\CC_l,\,l\in\QiG$, 
\begin{equation}\label{eq:neighbor_clique}
\bigcup_{l\in \QiG} \CC_l \subset \CN_i,
\end{equation}
holds \cite{sakurama2022generalized}.
Note that agent $i$ can independently obtain
$\CCl,\,l\in\QiG$ from the undirected subgraph induced by $\CN_i$.

\paragraph{Operator splitting}

Consider
\begin{align}\label{problem_3O}
\begin{array}{cl}
\underset{
\substack{y\in\BR^d}
}
{\mbox{min}}&
\displaystyle f(y) + g(y) + h(y),
\end{array}
\end{align}
where $f:\BR^{d}\to\BR$ is an smooth convex function, and
$g,h:\BR^{d}\to\BR\cup\{\infty\}$ are proper, closed, and convex functions.
For this problem, the following versatile algorithm, called \textit{(variable metric) Davis-Yin splitting} (DYS), has been proposed in \cite{davis2017three}:
\begin{align}\label{alg:DYS_variable_metric_pre}
\!\!\!\!\!\!\!\!\!\!\!\!\!
\begin{array}{lll}
    &y^{k+1/2} = \prox_{\alpha h}^M (z^k)\\
    &y^{k+1} =\prox_{\alpha g}^M(2 y^{k+1/2}- z^k - \alpha M^{-1} \nabla f(y^{k+1/2}) )\\
    & z^{k+1} = z^{k} + y^{k+1} - y^{k+1/2},
\end{array}
\end{align}
where $M\in\BR^{d\times d}$ is a positive definite symmetric matrix.
Note that the case of $M=I$ corresponds to the standard DYS.
This algorithm reduces to the Douglas-Rachford splitting when $f=0$ and to forward-backward splitting when $g=0$.
We have the following basic result, which states that $y^{k+1/2}$ and $y^{k+1}$ converge to a solution to \eqref{problem_3O} under an appropriate $\alpha>0$.
For further convergence results, see 
\cite{davis2017three,ryu2022large,aragon2022direct,lee2022convergence,yi2022convergence,condat2022randprox} and Subsection \ref{subsec:linear-rate}.
\begin{lemma}\label{lem:DYS}
Suppose that $M^{-1/2}\nabla f(y)M^{-1/2}$ is $L$-Lipschitz continuous for positive definite $M$.
Let $z^0\in\BR^d$ and $\alpha\in(0,2/L)$.
Assume that Problem \eqref{problem_3O} has an optimal solution.
Then, $y^k$ and $y^{k+1/2}$ updated by \eqref{alg:DYS_variable_metric_pre} converge to an optimal solution to Problem \eqref{problem_3O}.
\end{lemma}

\section{Clique-Wise Duplication Matrix}\label{sec:CDmatrix}

This section presents the definition and properties of the CD matrix $\BD$ that allows us to leverage operator splitting techniques for Problem \eqref{problem} in a distributed fashion\footnote{Note that the matrix $\BD$ itself is not new.
The same or similar ideas can be found in other existing papers, e.g., SDP \cite{vandenberghe2015chordal,zheng2019distributed} and generalized Nash equilibrium seeking \cite{bianchi2024end}.}
We also present a new mixing matrix $\BPhi$ with the matrix $\BD$, showing a preferable distribution of eigenvalues.

\subsection{Fundamentals}

The definition and essential properties of the CD matrix are presented in what follows.

First, we assume the non-emptiness of $\QiGactive$.
If this assumption is not satisfied, we can alternatively consider a subgraph induced by the node set to $\bigcup_{l\in\QGactive} \CC_l$.
\begin{assu}\label{assumption:sumCl=N}
For all $i\in\CN$, $\QiGactive\neq\emptyset$ holds.
\end{assu}

Then, the definition of the CD matrix is given as follows.
Here, $d_i$ for each $i\in\CN$ is the size of $x_i$ in Problem \eqref{problem}, and we define 
\begin{align*}
d = \sum_{i=1}^n d_i,
    \quad
    d^l = \sum_{j\in\CCl}d_j,\quad
    \hat{d}=\sum_{l\in\QGactive}d^l.
\end{align*}

\begin{defi}\label{def:CDmatrix}
For $d_i,\,i\in\CN$ and cliques $\CCl,\,l\in\QGactive$ of graph $\CG$, the \textit{Clique-wise Duplication (CD) matrix} $\BD$ is defined as
\begin{equation}
    \BD := 
    \begin{bmatrix}
        D_1\\
        \vdots\\
        D_{\qstar}
    \end{bmatrix}
    \in \BR^{\hat{d} \times d},
\end{equation}
where $D_l = [E_j]_{j\in\CCl}\in \BR^{d^l \times d}$ and $E_j = [O_{d_j\times d_1},\ldots,I_{d_j},\ldots,O_{d_j\times d_n}] \in \BR^{d_j\times d}$
for each $l\in\QGactive$.
\end{defi}

The CD matrix $\BD$ can be interpreted as follows. For $\Bx=[x_1^\top,\ldots,x_n^\top]^\top\in\BR^d$, $\BD \Bx = [x_\CCl]_{l\in\QGactive} \in \BR^{\hat{d}}$
holds since $D_l \Bx = x_\CCl \in \BR^{d^l}$. Hence, the CD matrix $\BD$ generates the copies of $\Bx$ with respect to cliques $\CCl,\,l\in\QGactive$.

The following lemma provides the fundamental properties of the CD matrix.
Now, let the matrix $E_{l,i}\in \BR^{d_i\times d^l}$ be 
\begin{equation}\label{Eli}
    E_{l,i} = [O_{d_i\times d_{j_1}},\ldots,I_{d_i},\ldots,O_{d_i\times d_{j_{|\CCl|}}}] \in \BR^{d_i\times d^l}
\end{equation}
for $\CCl=\{j_1,\ldots,i,\ldots,j_{|\CCl|}\},\,l\in\QiGactive$.
This matrix $E_{l,i}$ fulfills 
$E_{l,i} x_{\CCl} = x_i$
for $x_\CCl$ and $i\in\CCl$.
\begin{lemma}\label{thorem:CDmatrix}
    Under Assumption \ref{assumption:sumCl=N}, the followings hold.
    \begin{itemize}
        \item [(a)] $\BD$ is column full rank.
        \item[(b)] $\BD^\top\BD = \bdiag (|\CQ_\CG^1| I_{d_1},\ldots,|\CQ_\CG^n| I_{d_n}) \succ O$.
        \item[(c)] For $\By =[y_l]_{l\in\QGactive} \in \BR^{\hat{d}}$ with $y_l\in\BR^{d^l}$,
        \begin{equation}\label{eq:DTtimesy}
            \BD^\top \By = 
            \begin{bmatrix}
            \sum_{l\in \CQ_\CG^1} E_{l,1} y_l\\
            \vdots\\
            \sum_{l\in \CQ_\CG^n} E_{l,n} y_l
            \end{bmatrix}
            \in \BR^d.
        \end{equation}
    \end{itemize}
\end{lemma}

Using the CD matrix and \eqref{eq:neighbor_clique}, we can distributedly compute the least squares solution of $\By=\BD\Bx$, i.e.,
\begin{align}\label{eq:x_y}
    \Bx = (\BD^\top\BD)^{-1}  \BD^\top \By
\end{align}
and the projection of $\By$ onto $\mathrm{Im}(\BD)$ as
\begin{align}\label{eq:proj}
    P_{\mathrm{Im}(\BD)}(\By)= \BD(\BD^\top\BD)^{-1}  \BD^\top \By.
\end{align}

\begin{exam}\label{ex:CDmatrix}
Consider $\CG=(\CN,\CE)$ with $\CN=\{1,2,3\}$ and $\CE=\{\{1,2\},\,\{2,3\}\}$.
Let $d_1=d_2=d_3=1$ and $\QGactive = \{1,2\}$ with $\CC_1=\{1,2\}$ and $\CC_2=\{2,3\}$.
Then, we obtain $\QGactive^1=\{1\}$, $\QGactive^2 = \{1,2\}$, and $\QGactive^3=\{2\}$, 
which ensures Assumption \ref{assumption:sumCl=N}. 
For this system, the CD matrix is given by $\BD = [D_1^\top,D_2^\top]^\top \in \BR^{4 \times 3}$, where
$D_1 = 
\begin{bmatrix}
1 & 0 & 0\\
0& 1 & 0 
\end{bmatrix},\,
D_2 = 
\begin{bmatrix}
0 & 1 & 0\\
0& 0 & 1 
\end{bmatrix}$.
We then obtain $D_1\Bx= [x_1,x_2]^\top$ and $D_2\Bx = [x_2,x_3]^\top$ for $\Bx=[x_1,x_2,x_3]^\top\in\BR^3$. Moreover, 
$\BD^\top\BD = D_1^\top D_1 + D_2^\top D_2 =\diag(1,2,1) =\diag(|\QGactive^1|,|\QGactive^2|,|\QGactive^3|) $, and
\begin{align*}
\BD^\top \By = D_1^\top y_1 + D_2^\top y_2 = \begin{bmatrix}
    y_{1,1} \\
    y_{1,2} + y_{2,1}\\
    y_{2,2}
\end{bmatrix}
=
\begin{bmatrix}
    E_{1,1}y_{1} \\
    E_{1,2}y_{1} + E_{2,2} y_{2}\\
    E_{2,3} y_{2}
\end{bmatrix}
\end{align*}
for any vector $\By = [y_1^\top,y_2^\top]^\top\in\BR^4$ with $y_1=[y_{1,1},y_{1,2}]^\top\in\BR^2$ and $y_2=[y_{2,1},y_{2,2}]^\top\in\BR^2$,
which can be computed in a distributed fashion.
\end{exam}

\subsection{Useful properties}

Here, we provide useful properties of the CD matrix $\BD$.

The following result shows that the gradient and proximal operator with $\BD$ can be computed in a distributed fashion.
Here, $i$th block $x_i$ of $\Bx = (\BD^\top\BD)^{-1}\BD^\top \By$ is represented by
$x_i = E_i (\BD^\top\BD)^{-1}\BD^\top \By =  \frac{1}{|\QiGactive|} \sum_{l\in\QiGactive} E_{l,i}y_l$
from Lemma \ref{thorem:CDmatrix}.
\begin{prop}\label{thorem:CDmatrix_grad_prox}
    Let $\By\in \BR^{\hat{d}}$.
    Then, under Assumption \ref{assumption:sumCl=N}, the following equations hold.
    \begin{itemize}
        \item[(a)] Let $\hat{g}_i:\BR^{d_i}\to\BR\cup\{\infty\}$ be a proper, closed, and convex function for each $i\in\CN$.
        Define $G:\BR^{\hat{d}}\to\BR\cup\{\infty\}$ as
           $G(\Bz)=  \delta_{\mathrm{Im}(\BD)}(\Bz) + \sum_{i=1}^n \hat{g}_i(E_i (\BD^\top\BD)^{-1}\BD^\top \Bz) )$.
        Let $\alpha>0$. Then,
    \begin{align}\label{eq:prox_CDmatrix}
       \prox_{\alpha G}(\By)=\BD 
        \begin{bmatrix}
        \prox_{\frac{\alpha}{|\QGactive^1|}  \hat{g}_1} (
        E_1 (\BD^\top\BD)^{-1}\BD^\top \By)
        ) \\
        \vdots \\
         \prox_{\frac{\alpha}{|\QGactive^n|}  \hat{g}_n} (
         E_n (\BD^\top\BD)^{-1}\BD^\top \By)
         ) 
        \end{bmatrix}.
        \end{align}
        \item [(b)] 
        Let $\BQ=\bdiag([Q_l]_{l\in\QGactive})$, where $Q_l = \bdiag([\frac{1}{|\QjGactive|} I_{d_j}]_{j\in\CCl})$ for each $l\in\QGactive$. 
        Then,
    \begin{align}\label{eq:prox_CDmatrix_2}
       \prox_{\alpha G }^\BQ(\By)=\BD 
        \begin{bmatrix}
        \prox_{\alpha\hat{g}_1} (
        E_1 (\BD^\top\BD)^{-1}\BD^\top \By)
        ) \\
        \vdots \\
         \prox_{\alpha \hat{g}_n} (
         E_n (\BD^\top\BD)^{-1}\BD^\top \By)
         ) 
        \end{bmatrix}.
        \end{align}
        \item[(c)] Let $\hat{f}_i:\BR^{d_i}\to\BR$ be a differentiable function. Then,
        \begin{align}\label{eq:grad_CDmatrix}
          \frac{\pd }{\pd \By} \sum_{i=1}^n \hat{f}_i(E_i (\BD^\top\BD)^{-1}\BD^\top \By ) 
            = \BD (\BD^\top\BD)^{-1}
            \begin{bmatrix}
        \nabla_{x_1} \hat{f}_1 (E_1(\BD^\top\BD)^{-1}\BD^\top \By ) \\
        \vdots\\
         \nabla_{x_n} \hat{f}_n (E_n(\BD^\top\BD)^{-1}\BD^\top \By ) \\
            \end{bmatrix}.
        \end{align}
    \end{itemize}
\end{prop}

Additionally, we provide properties of the CD matrix concerning matrix $\BQ$. 
Those properties are useful to derive the NIDS \cite{li2019decentralized} and Exact diffusion \cite{yuan2018exact} from the proposed method.  
\begin{prop}\label{prop:CDmatrix_variable_metric}
Let $\BQ$ denote the matrix in Proposition \ref{thorem:CDmatrix_grad_prox}b.
Then, under Assumption \ref{assumption:sumCl=N}, the following equations hold:
\begin{itemize}
    \item[(a)] $\BD^\top\BQ\BD = I_d$.
    \item[(b)] $\BD^\top \BQ = (\BD^\top\BD)^{-1}\BD^\top$ and $\BD^\top \BQ^{-1} = \BD^\top\BD\BD^\top$.
    \item[(c)] $\BQ\BD = \BD (\BD^\top\BD)^{-1}$ and $\BQ^{-1}\BD = \BD \BD^\top\BD$.
\end{itemize}   
\end{prop}

\subsection{A mixing matrix $\BPhi$}\label{subsec:Phi}

Using the CD matrix and the matrices $Q_l,\,l\in\QGactive$ in Proposition \ref{thorem:CDmatrix_grad_prox}b, we can obtain a positive semidefinite and doubly stochastic matrix $\BPhi$ that will be used in Section \ref{sec:analysis}.
Thanks to the definition, $\BPhi$ in \eqref{eq:clique_doubly_stochastic} can be constructed only from local information (i.e., $\QjGactive,\,j\in\bigcup_{l\in\QiGactive}\CCl\subset\CN_i$).
Note that this matrix can be viewed as a special case of the clique-based projection $T$ in \cite{watanabe2022distributed} and Appendix \ref{A:CPGD} for the consensus constraint, i.e., $\BPhi \Bx= \BD^\top (\BD^\top\BD)^{-1}P^\BQ_{\Pi_{l\in\QGactive}\CD_l}(\BD\Bx)$ for $\CD_l$ in \eqref{D_l}.
We here pose the following assumption\footnote{Assumption \ref{assumption:edge} is not strict and satisfied for $\QGactive=\QG$, $\maxQG$, $\QGactive^\mathrm{edge}$.}.
\begin{assu}\label{assumption:edge}
For $\QGactive$,
$\QiGactive\cap \QjGactive \neq \emptyset \Leftrightarrow \{i,j\}\in\CE$.
\end{assu}

The matrix $\BPhi$ and its basic properties are given as follows.
\begin{prop}\label{prop:double_stochastic}
Suppose Assumptions \ref{assumption:sumCl=N} and \ref{assumption:edge}.
Consider the matrices $Q_l,\,l\in\QGactive$ in Proposition \ref{thorem:CDmatrix_grad_prox}b.
Suppose that $d_1=\cdots=d_n=1$.
Then, 
\begin{equation}\label{eq:clique_doubly_stochastic}
\BPhi =
\begin{bmatrix}
   \frac{1}{|\QGactive^1|}\sum_{l\in\QGactive^1}
\frac{\mathbf{1}_{|\CCl|}^\top Q_l D_l}{\mathbf{1}_{|\CCl|}^\top Q_l \mathbf{1}_{|\CCl|}}\\
   \vdots\\
   \frac{1}{|\QGactive^n|}\sum_{l\in\QGactive^n}
\frac{\mathbf{1}_{|\CCl|}^\top Q_l D_l}{\mathbf{1}_{|\CCl|}^\top Q_l \mathbf{1}_{|\CCl|}}\\
\end{bmatrix}
\in \BR^{n\times n}
\end{equation}
is doubly stochastic,  
and it holds that
\begin{equation}\label{eq:Phi_ij}
[\BPhi]_{ij} =
\begin{cases}
    \frac{1}{|\QiGactive||\QjGactive|} \sum_{l\in\QGactive^{ij}}  \frac{1}{\mathbf{1}_{|\CCl|}^\top Q_l \mathbf{1}_{|\CCl|}}, & \{i,j\}\in \CE
    \\
    0, & \text{otherwise},
\end{cases}
\end{equation}
where $[\BPhi]_{ij}$ represents $(i,j)$ entry of $\BPhi$.
Moreover, $\lambda_\mathrm{max}(\BPhi)=1$ and $\lambda_\mathrm{min}(\BPhi)\geq 0$ hold.
Furthermore, when $\CG$ is connected, the eigenvalue $1$ of $\BPhi$ is simple.
\end{prop}
\begin{proof}
The right stochasticity is proved as $(\frac{1}{|\QGactive^i|}\sum_{l\in\QGactive^i}
\frac{\mathbf{1}_{|\CCl|}^\top Q_l D_l}{\mathbf{1}_{|\CCl|}^\top Q_l \mathbf{1}_{|\CCl|}} )\mathbf{1}_n 
=      \frac{1}{|\QGactive^i|}\sum_{l\in\QGactive^i}
\frac{\mathbf{1}_{|\CCl|}^\top Q_l\mathbf{1}_{|\CCl|}}{\mathbf{1}_{|\CCl|}^\top Q_l \mathbf{1}_{|\CCl|}}
= \frac{1}{|\QGactive^i|}\sum_{l\in\QGactive^i} 1 = 1.$
Using the definition of $D_l$ in Definition \ref{def:CDmatrix}, the left stochasticity is also verified as
\begin{align*}
&
\mathbf{1}_n^\top \BPhi
=\sum_{i=1}^n\frac{1}{|\QGactive^i|}\sum_{l\in\QGactive^i}
\frac{\mathbf{1}_{|\CCl|}^\top Q_l D_l}{\mathbf{1}_{|\CCl|}^\top Q_l \mathbf{1}_{|\CCl|}}
\\
=
&
\sum_{l\in\QGactive}
\sum_{j\in\CCl} \frac{1}{|\QGactive^j|}
\frac{\mathbf{1}_{|\CCl|}^\top Q_l D_l}{\mathbf{1}_{|\CCl|}^\top Q_l \mathbf{1}_{|\CCl|}}
=
\sum_{l\in\QGactive}\sum_{j\in\CCl} \frac{1}{|\QjGactive|} E_j
\\
=
&
\sum_{i=1}^n \frac{1}{|\QiGactive|} \sum_{l\in\QiGactive} E_i = \sum_{i=1}^n E_i = \mathbf{1}_n^\top
\end{align*}
from $\mathbf{1}_{|\CCl|}^\top Q_l \mathbf{1}_{|\CCl|}$ and
$\mathbf{1}_{|\CCl|}^\top Q_l D_l = \sum_{j\in\CCl} \frac{1}{|\QjGactive|}E_j$.
Next, 
\begin{align*}
[\BPhi]_{ij} 
= E_i \BPhi E_j^\top= \frac{1}{|\QGactive^i|}\sum_{l\in\QGactive^i}
\frac{\mathbf{1}_{|\CCl|}^\top Q_l D_l}{\mathbf{1}_{|\CCl|}^\top Q_l \mathbf{1}_{|\CCl|}}
E_j^\top 
= \frac{1}{|\QGactive^i|}\sum_{l\in\QGactive^i}
\sum_{p\in\CC_l}\frac{1}{\mathbf{1}_{|\CCl|}^\top Q_l \mathbf{1}_{|\CCl|}} \frac{1}{|\QGactive^p|}E_p
E_j^\top
\end{align*}
holds. Then, we obtain \eqref{eq:Phi_ij}.
Moreover, $\lambda_\mathrm{max}(\BPhi)=1$ directly follows from Gershgorin disks theorem \cite{bullo2020lectures}. Additionally, from the firmly nonexpansiveness of the clique-based projection $T$ (see Proposition \ref{prop:clique-based_projection}), we obtain $\Bx^\top \BPhi \Bx \geq \|\BPhi\Bx\|^2$ for any $\Bx\in\BR^n$, which gives $\lambda_\mathrm{min}(\BPhi)\geq 0$.

Finally, by the assumption of $\QGactive^{ij} \neq \emptyset \Leftrightarrow \{i,j\}\in\CE$, the associated graph of $\BPhi$ is equal to $\CG$.
Therefore, the eigenvalue $1$ of $\BPhi$ is simple
when $\CG$ is connected (see \cite{bullo2020lectures}).
\end{proof}

Now, we state the following proposition for $\BPhi$ in \eqref{eq:clique_doubly_stochastic}, implying that $\BPhi$ enables fast consensus.
This is because all the eigenvalues smaller than $1$ are likely to take smaller values than other popular positive semidefinite mixing matrices from the Gershgorin disks theorem \cite{bullo2020lectures}.
A numerical example of the eigenvalues and a sketch of Proposition \ref{prop:Phi_eigval}'s implication are illustrated in Fig. \ref{fig:W_comparison}.
\begin{prop}\label{prop:Phi_eigval}
Suppose Assumption \ref{assumption:sumCl=N}.
For undirected connected graph $\CG$,
consider the matrix $\BPhi$ in \eqref{eq:clique_doubly_stochastic} with $\QGactive=\QGactive^{\mathrm{edge}}$.
Let $\widetilde{\BW}_{\BL}=(I+\BW_\BL)/2$, where $\BW_\BL = I-\epsilon \BL$ with
Laplacian matrix $\BL$ of $\CG$ and $\epsilon \in(0,1/\max_{i\in\CN}\{|\CN_i|-1\} )$\footnote{
The matrix $\BL$ is defined as $[\BL]_{ii}=|\CN_i|-1$ for $i=1,\ldots,n$ and $[\BL]_{ij} = -1$ for $\{i,j\}\in\CE$. Otherwise $[\BL]_{ij}=0$.
In \cite{xiao2006distributed}, $[\BW_{\BL}]_{ij}$ with $\epsilon=1/\max_{i\in\CN}\{|\CN_i|\}$ is said to be \textit{the max-degree weight}.
}.
Similarly, let  $\widetilde{\BW}_{\mathrm{mh}}=(I+\BW_\mathrm{mh})/2$ with $\BW_\mathrm{mh}$ obtained by the Metropolis–Hastings weights\footnote{$\BW_\mathrm{mh}$ is defined as $[\BW_\mathrm{mh}]_{ij}=1/(\max\{|\CN_i|-1,|\CN_j|-1\}+\varepsilon)$ for $\{i,j\}\in\CE$ and $[\BW_\mathrm{mh}]_{ii}=1-\sum_{j\in\CN_i\setminus\{i\}}[\BW_\mathrm{mh}]_{ij}$ for $i=1,\ldots,n$. Otherwise $[\BW_\mathrm{mh}]_{ij}=0$.} in \cite{xiao2006distributed}.
Then, for all $i=1,\ldots,n$, we have \begin{equation}\label{Phi_diagonal}
    [\BPhi]_{ii} < [\widetilde{\BW}_{\BL}]_{ii},\quad
    [\BPhi]_{ii} < [\widetilde{\BW}_{\mathrm{mh}}]_{ii}.
\end{equation}
\end{prop}
 \begin{proof}
When $\QGactive=\QGactive^{\mathrm{edge}}$, $|\QiGactive| = |\CN_i|-1$ holds for $i=1,\ldots,n$, and $\QGactive^{ij}$ for $\{i,j\}\in\CE$ becomes a singleton $\QGactive^{ij}=\{l\}$ with $\CCl=\{i,j\}$ 
as $\QiGactive$ is just the set of indices of edges that include $i$.
Then, for $\{i,j\}\in\CE$ we get $[\BPhi]_{ij} = \frac{1}{|\CN_i|-1+|\CN_j|-1}$.
Hence, recalling the definition of $\widetilde{\BW}_\BL$ and $\widetilde{\BW}_\mathrm{mh}$ for $\{i,j\}\in\CE$, we have
$[\widetilde{\BW}_\BL]_{ij}= 1/2\epsilon < 1/(2\max_{i\in\CN}\{|\CN_i|-1\})\leq[\BPhi]_{ij}$ and $[\widetilde{\BW}_\mathrm{mh}]_{ij} = 1/2(\max\{|\CN_i|-1,|\CN_j|-1\}+\varepsilon)<1/2\max\{|\CN_i|-1,|\CN_j|-1\}\leq [\BPhi]_{ij}$, respectively.
Therefore, since all $(i,j)$ entries of $\BPhi$ for $\{i,j\} \in\CE$ are bigger than those of $\widetilde{\BW}_\BL$ and $\widetilde{\BW}_\mathrm{mh}$ and these matrices are doubly stochastic, we get \eqref{Phi_diagonal}.
 \end{proof}
\begin{rem}
In Fig. \ref{fig:W_comparison}a, we use not $\QGactive^\mathrm{edge}$ but $\maxQG$, which also realizes smaller eigenvalues.
Likewise,
even when $\QGactive \neq \QGactive^{\mathrm{edge}}$, $\BPhi$ can have smaller eigenvalues than $\widetilde{\BW}_\BL$ and $\widetilde{\BW}_\mathrm{mh}$.
\begin{figure}
    \centering
    \includegraphics[width=0.8\columnwidth]{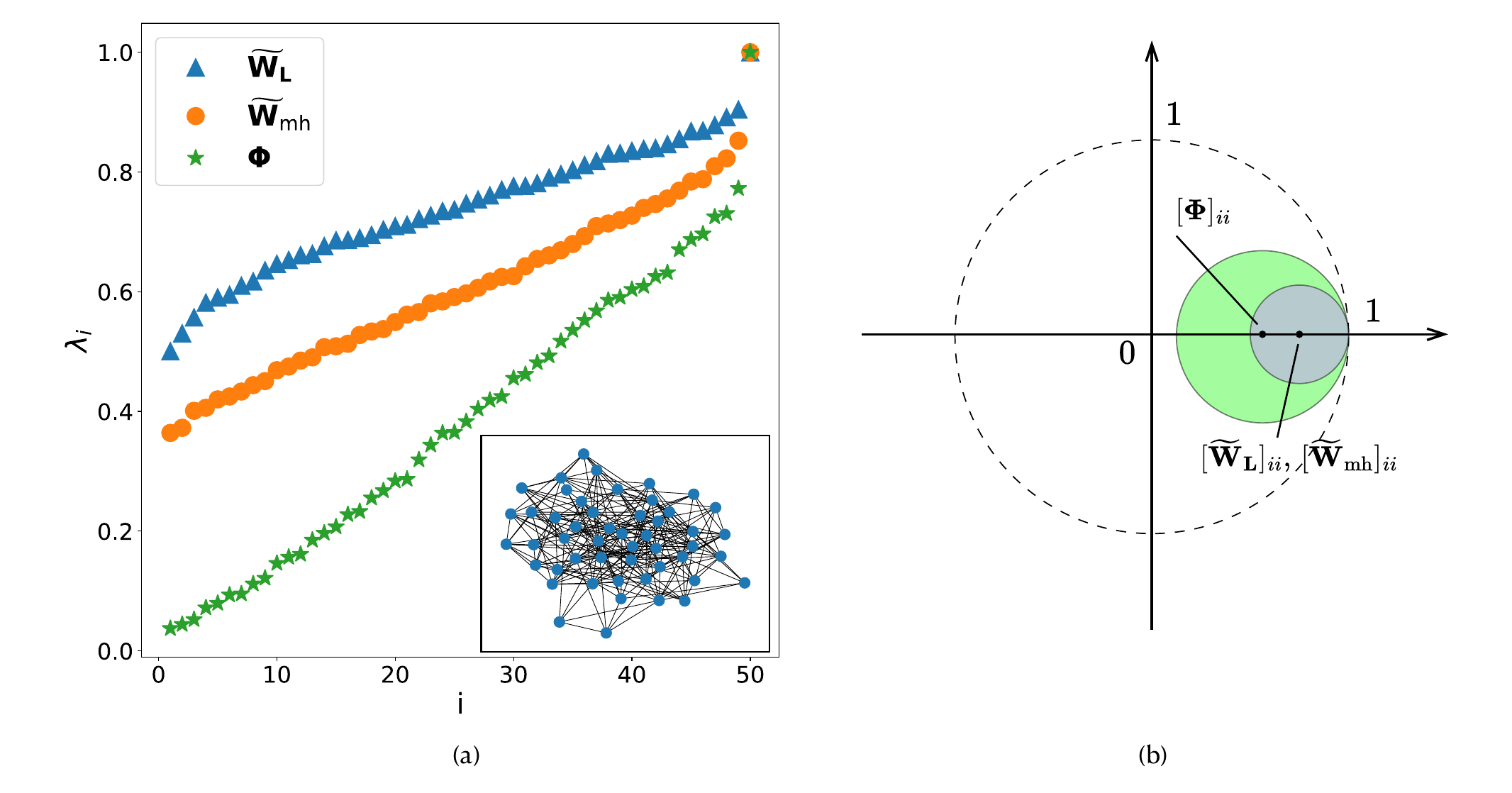}
    \caption{(a) Comparison of the eigenvalues $\lambda_i$ of $\BPhi$, $\widetilde{\BW}_\BL$ with $\epsilon=0.99/ (\max_{i\in\CN}|\CN|_i-1)$, and $\widetilde{\BW}_\mathrm{mh}$
    for a random graph with $n=50$ inside the plot; (b) A sketch of the region of each eigenvalue of $\BPhi$, $\widetilde{\BW}_\BL$, and $\widetilde{\BW}_\mathrm{mh}$ that Proposition \ref{prop:Phi_eigval} and the Gershgorin disks theorem imply.
    Both indicate that $\BPhi$ probably takes smaller eigenvalues.}
    \label{fig:W_comparison}
\end{figure}
\end{rem}

\section{
Solution to Clique-Wise Coupled Problems via Operator Splitting
}\label{sec:alg}

This section presents our proposed algorithms for Problem \eqref{problem} with the CD matrix and DYS in \eqref{alg:DYS_variable_metric_pre} with the metrics of $M=I$ and $M=\BQ$.
We now assume the following.
\begin{assu}\label{assu:prob}
Problem \eqref{problem} has an optimal solution.
\end{assu}
In what follows, the functions $f:\BR^{\hat{d}}\to\BR$, $g:\BR^{\hat{d}}\to\BR$, $\hat{f}:\BR^{d}\to\BR$, and $\hat{g}:\BR^{d}\to\BR$ represent 
\begin{align}\label{f_g}
    & f(\By) = \sum_{l\in\QGactive}f_l(y_l),\quad g(\By) = \sum_{l\in\QGactive}g_l(y_l),\\
    & \hat{f}(\Bx) = \sum_{i=1}^n f_i(x_i),\quad
    \hat{g}(\Bx) = \sum_{i=1}^n g_i(x_i).
 \end{align}

\subsection{Algorithm description}\label{subsec:development}

\begin{algorithm}[t]
\caption{
Clique-based distributed Davis-Yin splitting (CD-DYS) algorithm for agent $i\in\CN$.
}
\label{alg:DYS_distributed_aw}
\begin{algorithmic}[1]
\REQUIRE $z_l^0$ and $\alpha>0$ for all $l\in\QiGactive$.
\FOR{$k = 0,1,\ldots$}
\STATE $x_i^k = \prox_{\frac{\alpha}{|\QGactive^i|}  \hat{g}_i} (
        \frac{1}{|\QiGactive|} \sum_{l\in\QiGactive} E_{l,i}z_l^k) $
\STATE \textcolor{black}{Gather $x_j^k$ from the neighboring agents $j\in\bigcup_{l\in\QiGactive}\CCl \subset \CN_i$.}
\STATE Obtain $y_l^{k+1/2}$, $y_l^{k+1}$, and $z_l^{k+1}$ for $l\in\QiGactive$ by
\begin{align*}
   &y_l^{k+1/2} = x^k_\CCl\\
    &y_l^{k+1} = \prox_{\alpha g_l} (2y_l^{k+1/2}-z_l^k - \alpha \nabla_{y_l} f_l(y_l^{k+1/2}) 
    - \alpha 
    [\frac{1}{|\QGactive^j|} \nabla_{x_j} \hat{f}_j(x_j^k)]_{j\in\CCl}
    )
    \\
    &z_l^{k+1} = z_l^k + y_l^{k+1} - y_l^{k+1/2}
\end{align*}
\ENDFOR
\end{algorithmic}
\end{algorithm}

We give the distributed optimization algorithm in Algorithm \ref{alg:DYS_distributed_aw}, the \textit{clique-based distributed Davis-Yin splitting (CD-DYS)} algorithm.
This algorithm is distributed from \eqref{eq:neighbor_clique}.
By Lemma \ref{thorem:CDmatrix},
the aggregated form of this algorithm is as follows:
\begin{align}\label{alg:DYS}
\!\!\!\!\!\!\!\!\!
\begin{array}{lll}
    &\Bx^k = \prox_{\sum_{i=1}^n \frac{\alpha}{|\QiGactive|}\hat{g}_i(\cdot)}(\BD^\top\BD)^{-1}\BD^\top \Bz^k \\
    &\By^{k+1/2} = \BD \Bx^k\\
    &\By^{k+1} = \prox_{\alpha g} (2\By^{k+1/2}-\Bz^k  
    - \alpha \nabla_{\By} f (\By^{k+1/2})- \alpha\BD(\BD^\top\BD)^{-1}\nabla_\Bx \hat{f}(\Bx^k)
    )
    \\
    &\Bz^{k+1} = \Bz^k + \By^{k+1} - \By^{k+1/2},
\end{array}
\end{align}
where $\Bx^k=[x_1^{k\top},\ldots,x_n^{k\top}]^\top $,
$\By^k= [y_l^k]_{l\in\QGactive}$, $\By^{k+1/2}= [y_l^{k+1/2}]_{l\in\QGactive}$, and $\Bz^k= [z_l^k]_{l\in\QGactive}$.
By Lemma \ref{lemma:convergence}, this algorithm converges to the optimal point under $\alpha\in(0,2/({\max_{l\in\QGactive}L_l+\max_{i\in\CN}\frac{\hat{L}_i}{|\QiGactive|}}))$.

This algorithm can be derived in the following way.
From \eqref{eq:x_y},
for $\By=\BD\Bx \in\mathrm{Im}(\BD)$,
we can reformulate Problem \eqref{problem} into the form of \eqref{problem_3O} as follows:
\begin{align}\label{problem_y_aw}
\begin{array}{cl}
&\underset{
\substack{y_l\in \BR^{d^l},\,l\in\QGactive}
}
{\mbox{min}}  \quad
\displaystyle  \underbrace{\sum_{i=1}^n \hat{f}_i (E_i (\BD^\top\BD)^{-1}\BD^\top \By) +\sum_{l\in\QGactive} f_l (y_l)}_{f \text{ in \eqref{problem_3O}}} 
\\
&\displaystyle  
+ \underbrace{\sum_{l\in\QGactive} g_l (y_l) }_{g \text{ in \eqref{problem_3O}}}
\displaystyle
+ \underbrace{\sum_{i=1}^n \hat{g}_i (E_i (\BD^\top\BD)^{-1}\BD^\top \By)  + \delta_{\mathrm{Im}(\BD)}(\By)}_{h \text{ in \eqref{problem_3O}}}.
\end{array}
\!\!\!\!
\end{align}
For Problem \eqref{problem_y_aw}, Proposition \ref{thorem:CDmatrix_grad_prox} tells us that
the function $\sum_{i=1}^n \hat{g}_i (E_i (\BD^\top\BD)^{-1}\BD^\top \By)  + \delta_{\mathrm{Im}(\BD)}(\By)$ is proximable for proximable $\hat{g}_i$, and the proximal operator can be computed in a distributed fashion.
Accordingly, we can directly apply DYS in \eqref{alg:DYS_variable_metric_pre} with $M=I$ to \eqref{problem_y_aw}. 
From Proposition \ref{thorem:CDmatrix_grad_prox},
setting $x_i^k=
 \prox_{\frac{\alpha}{|\QGactive^i|}  \hat{g}_i} (
        E_i (\BD^\top\BD)^{-1}\BD^\top \Bz^k) $
gives the distributed algorithm in Algorithm \ref{alg:DYS_distributed_aw} (or \eqref{alg:DYS}).
To implement Algorithm \ref{alg:DYS_distributed_aw}, the gradient information $\nabla_{x_j} \hat{f}_j$
has to be shared in the neighbors.
Provided the agents are collaborative within their neighbors, this requirement is not restrictive.
We also note that in special cases as consensus optimization, this requirement can be alleviated; see Section \ref{sec:analysis}.

\subsection{Variable metric version}
Applying the variable metric DYS in \eqref{alg:DYS_variable_metric_pre} with $M=\BQ$ in Proposition \ref{thorem:CDmatrix_grad_prox} instead to Problem \eqref{problem_y_aw} gives the following algorithm:
\begin{align}\label{alg:DYS_distributed_aw_vm}
\!\!\!\!\!\!\!\!\!\!\!\!\!\!
\begin{array}{lll}
&x_i^k = \prox_{\alpha \hat{g}_i} (
        \frac{1}{|\QiGactive|} \sum_{l\in\QiGactive} E_{l,i}z_l^k)\\
   &y_l^{k+1/2} = x^k_\CCl\\
    &y_l^{k+1} = \prox_{\alpha g_l}^{Q_l} (2y_l^{k+1/2}-z_l^k - \alpha Q_l^{-1} \nabla_{y_l} f_l(y_l^{k+1/2}) 
    - \alpha 
    [\nabla_{x_j} \hat{f}_j(x_j^k)]_{j\in\CCl})
    \\
    &z_l^{k+1} = z_l^k + y_l^{k+1} - y_l^{k+1/2},
\end{array}
\end{align}
where we have used Propositions \ref{thorem:CDmatrix_grad_prox}b and \ref{prop:CDmatrix_variable_metric}.
It will turn out in Section \ref{sec:analysis}
that this algorithm shows an interesting connection to other methods as Fig. \ref{fig:relation} through $\BPhi$ in \eqref{eq:clique_doubly_stochastic}. 
By Lemma \ref{lemma:convergence}, a sufficient condition for the convergence is $\alpha\in(0,2/({\max_{l\in\QGactive}\max_{j\in\CCl}(|\QjGactive|L_l)+\max_{i\in\CN}\hat{L}_i}))$.

\section{Revisit of Consensus Optimization as A Clique-Wise Coupled Problem}\label{sec:analysis}

\begin{figure}[t]
\centering
\begin{tikzpicture}
[
every node/.style={outer sep=0.12cm, inner sep=0},
arrow/.style={-{Stealth[length=0.25cm]}, thick},
block/.style={rectangle, draw, minimum height = 1.0cm,
minimum width=1.4cm, thick, outer sep = 0},
sum/.style={thick, circle, draw, inner sep=0,
minimum size=6pt, outer sep=0},
point/.style={radius=2pt}
]
\node [block,align=center] (a) at (0,0) {
Variable metric CD-DYS\\
(Alg. \eqref{alg:DYS_variable_metric_pre} for \eqref{problem} via \eqref{problem_y_aw})
};
\node [block,align=center] (b) at (-5,-2) {
Variable metric CD-DYS\\
w.r.t. $\BQ$ (Alg. \eqref{alg:DYS_distributed_aw_vm})};
\node [block] (c) at (5,-2) {
CD-DYS (Alg. \ref{alg:DYS_distributed_aw} or \eqref{alg:DYS}
};
\node [block] (d) at (-5,-4) {
CPGD \cite{watanabe2022distributed}};
\node [block] (e) at (5,-4) {
NIDS \cite{li2019decentralized} (Alg. \eqref{NIDS_prox} with $\widetilde{\BW}=\BPhi$) };
\node [block,align=center] (f) at (5,-6) {
Exact diffusion \cite{yuan2018exact} with $\widetilde{\BW}=\BPhi$};
\node [block] (g) at (-5,-6) {
Diffusion \cite{sayed2014diffusion} with $\widetilde{\BW}=\BPhi$};

\draw[arrow] (a) -- (b) node [left, pos=0.4] {\!\!\!\!\!\!\!\!\!$M=\BQ$};
\draw[arrow] (a) -- (c) node [right, pos=0.4] {\quad\quad$M=I$};
\draw[arrow] (b) -- (d) node [left, pos=0.5,align=center] {$g_l=\delta_{\CD_l},\,\hat{g}_i=0$ \\
+ approximation
};
\draw[arrow] (b) -- (e) node [right, pos=0.4,align=center] {\quad\quad$g_l=\delta_{\CD_l}$ 
with $\CD_l=$\eqref{D_l}
};
\draw[arrow] (e) -- (f) node [right, pos=0.4,align=center] {$\hat{g}_i=0$};
\draw[arrow] (f) -- (g) node [above, pos=0.5,align=center] {Approximation};
\draw[arrow] (d) -- (g) node [left, pos=0.5,align=center] {$\CD_l=$ \eqref{D_l}};
\end{tikzpicture}
\caption{The relationships between the proposed methods and existing methods.
}
\label{fig:relation}
\end{figure}

This section is dedicated to a detailed analysis of the proposed methods in Section \ref{sec:alg} in light of consensus optimization, presenting a renewed perspective.
We here demonstrate the relationship summarized in Fig. \ref{fig:relation}
by showing the most important part: Algorithm \eqref{alg:DYS_distributed_aw_vm}
generalizes the NIDS in \cite{li2019decentralized}.
Our analysis reveals that matrix $\BPhi$ in \eqref{eq:clique_doubly_stochastic} naturally arises in the NIDS.
This fact and Proposition \ref{prop:Phi_eigval} imply that the proposed algorithm enables fast convergence \cite{li2019decentralized} beyond standard mixing matrices.
Furthermore, we present a new linear convergence rate for the NIDS with a non-smooth term based on its DYS structure. The linear rate follows from a more general result for the DYS \eqref{alg:DYS_variable_metric_pre}.

We here consider a special case of Problem \eqref{problem} given by
\begin{align}\label{problem_clique-wise_constraints}
\!\!\!\!\!\!\!
\begin{array}{cl}
\underset{
\substack{x_i\in \BR^{d_i},\,i\in\CN}
}
{\mbox{min}}&
\displaystyle  \sum_{i=1}^n \hat{f}_i (x_i)+ \sum_{i=1}^n \hat{g}_i (x_i) +\sum_{l\in\QGactive} \underbrace{\delta_{\CD_l}(x_\CCl)}_{g_l(x_\CCl)}.
\end{array}
\!\!\!
\end{align}
When $m=d_1=\cdots=d_n$ and 
\begin{equation}\label{D_l}
    \CD_l = \{x_\CCl\in \BR^{|\CCl|m}: \exists \theta \in \BR^m \text{ s.t. } x_\CCl = \mathbf{1}_{|\CCl|} \otimes \theta \},
\end{equation}
this problem is called a \textit{consensus optimization problem}, which we discuss here.
Notice that according to \cite{sakurama2022generalized},
$\cap_{l\in\QGactive}\{\Bx\in\BR^{nm}:x_\CCl\in D_l\} = \{\Bx\in\BR^{nm}: x_1=\cdots=x_n\} $ holds under the connectivity of graph $\CG$ and Assumptions \ref{assumption:sumCl=N} and \ref{assumption:edge}.

Note also that the full analysis of Fig. \ref{fig:relation} is found in Appendix \ref{appendix:relation}.

\subsection{CD-DYS as generalized NIDS}\label{subsec:analysis_consensus}
\paragraph{NIDS algorithm}
First, the NIDS algorithm \cite{li2019decentralized} for consensus optimization for $k=1,2,\ldots$ is as follows:
\begin{align}\label{NIDS_prox}
\begin{array}{ll}
&\Bx^{k}= \prox_{\alpha \hat{g}} (\Bw^{k})\\
    & \Bw^{k+1}= \Bw^k - \Bx^{k} + \widetilde{\BW} (2\Bx^{k}-\Bx^{k-1} 
     +  \alpha \nabla_\Bx \hat{f}(\Bx^{k-1}) - \alpha \nabla_\Bx \hat{f}(\Bx^{k}))
\end{array}
\end{align}
with arbitrary $\Bx^0$ and $\Bw^{1} =  \Bx^0-\alpha \nabla_\Bx\hat{f}(\Bx^0)$.
The matrix $\widetilde{\BW}$ is a positive semidefinite doubly stochastic mixing matrix.
A standard choice of $\widetilde{\BW}$ with no use of global information is $\widetilde{\BW} = \widetilde{\BW}_\mathrm{mh}$ in Proposition \ref{prop:Phi_eigval}.
To make $\widetilde{\BW}$ less conservative, \cite{li2019decentralized} suggests that
some global information is necessary (e.g., the value of $\lambda_\mathrm{max}(\BW_\mathrm{mh})$).

\paragraph{Analysis}

We here present the following proposition, stating that the proposed algorithm in \eqref{alg:DYS_distributed_aw_vm} yields
the NIDS algorithm with $\widetilde{\BW}=\BPhi$ in \eqref{eq:clique_doubly_stochastic}, which can achieve fast convergence as shown in Proposition \ref{prop:Phi_eigval} and Fig. \ref{fig:W_comparison}.
Note that we show the case of $m=1$ for simplicity.
\begin{prop}\label{prop:NIDS}
Consider Algorithm \eqref{alg:DYS_distributed_aw_vm} for Problem \eqref{problem_clique-wise_constraints}.
Suppose Assumptions \ref{assumption:sumCl=N}--\ref{assu:prob}.
Assume that for all $i\in\CN$,
$\hat{f}_i:\BR^{d_i}\to\BR$ is $\hat{L}_i$-smooth and convex, and $\hat{g}_i:\BR^{{d_i}}\to\BR$
is proper, closed, and convex.
For arbitrary $\Bx^0$, let $\Bz^1=\BD(\Bx^0-\alpha \nabla_\Bx\hat{f}(\Bx^0))$.
Then, for $k=1,2,\ldots$,
$\Bw^{k} := (\BD^\top\BD)^{-1}\BD^\top\Bz^{k}$ and $\Bx^k := \prox_{\alpha\hat{g}}(\Bw^k)$ satisfy the update of NIDS in \eqref{NIDS_prox} with $\widetilde{\BW}=\BPhi$ in \eqref{eq:clique_doubly_stochastic}.
\end{prop}
\begin{proof}
By Lemma \ref{thorem:CDmatrix}b--c,
multiplying the third line of \eqref{alg:DYS_distributed_aw_vm} by $(\BD^\top\BD)^{-1}\BD^\top$ 
gives $w_i^{k+1} = w_i^{k} - x_i^k + \frac{1}{|\QiGactive|}\sum_{l\in\QiGactive}
    E_{l,i}\prox_{\delta_{\CD_l}}^{Q_l} (2x^{k}_\CCl-z_l^k - \alpha D_l \nabla_\Bx \hat{f}(\Bx^k))$.
Then, plugging in
$\prox_{\delta_{\CD_l}}^{Q_l}(x_\CCl) = P_{\CD_l}^{Q_l}(x_\CCl)= 
 \mathbf{1}_{|\CCl|} 
\frac{\mathbf{1}_{|\CCl|}^\top Q_l x_\CCl}{\mathbf{1}_{|\CCl|}^\top Q_l \mathbf{1}_{|\CCl|}}$,
\begin{align*}
    &w_i^{k+1} 
    = w_i^{k} - x_i^k
     + \frac{1}{|\QiGactive|}\sum_{l\in\QiGactive}
\frac{\mathbf{1}_{|\CCl|}^\top Q_l }{\mathbf{1}_{|\CCl|}^\top Q_l \mathbf{1}_{|\CCl|}}  (2x^{k}_\CCl 
    -z_l^k - \alpha  D_l \nabla_\Bx \hat{f}(\Bx^k) ). 
\end{align*}
Additionally, we can transform $ \mathbf{1}_{|\CCl|}^\top Q_l z_l^{k+1}$ for $k\geq 1$ into
\begin{align}\label{eq:z_v}
  \mathbf{1}_{|\CCl|}^\top Q_l z_l^{k+1} 
=& \mathbf{1}_{|\CCl|}^\top Q_l (z_l^k - x_\CCl^k) \nonumber +\mathbf{1}_{|\CCl|}^\top Q_l (2x^{k}_\CCl-z_l^k - \alpha  D_l \nabla_\Bx \hat{f}(\Bx^k) ) \nonumber\\
=&  \mathbf{1}_{|\CCl|}^\top Q_l (x_\CCl^k - \alpha  D_l\nabla_\Bx \hat{f}(\Bx^k) ).
\end{align}
Thus, for $k=1,2,\ldots$, recalling the initialization of $\Bz^1 = \BD(\Bx^0-\alpha\nabla_\Bx \hat{f}(\Bx^0))$ and
applying \eqref{eq:z_v} to $w_i^{k+1}$ provide
    $w_i^{k+1} 
    = w_i^{k} - x_i^k
    + \frac{1}{|\QiGactive|}\sum_{l\in\QiGactive}
\frac{\mathbf{1}_{|\CCl|}^\top Q_l }{\mathbf{1}_{|\CCl|}^\top Q_l \mathbf{1}_{|\CCl|}}
    (2x^{k}_\CCl
    -x^{k-1}_\CCl 
    + \alpha D_l  (\nabla_\Bx \hat{f}(\Bx^{k-1}) - \nabla_\Bx \hat{f}(\Bx^k)))
    = w_i^{k} - x_i^k 
     + \frac{1}{|\QiGactive|}\sum_{l\in\QiGactive}
\frac{\mathbf{1}_{|\CCl|}^\top Q_l D_l}{\mathbf{1}_{|\CCl|}^\top Q_l \mathbf{1}_{|\CCl|}}
    (2\Bx^{k}
    -\Bx^{k-1} 
    + \alpha (\nabla_\Bx \hat{f}(\Bx^{k-1}) - \nabla_\Bx \hat{f}(\Bx^k)))$.
Thus,
setting $\widetilde{\BW}=\BPhi$ with $\BPhi$ in \eqref{eq:clique_doubly_stochastic}, 
we get
$\Bw^{k+1} = \Bw^{k} - \Bx^k + \widetilde{\BW}(2\Bx^{k}
    -\Bx^{k-1}
    + \alpha (\nabla_\Bx \hat{f}(\Bx^{k-1}) -\alpha  \nabla_\Bx \hat{f}(\Bx^k)))$.
\end{proof}
\begin{rem}
The original NIDS paper \cite{li2019decentralized} states that the NIDS is obtained from the PD3O in \cite{yan2018new} (a primal-dual variant of DYS).
Meanwhile, in Proposition \ref{prop:NIDS}, we rely only on the primal part and obtain $\widetilde{\BW}=\BPhi$ as a fixed parameter.
\end{rem}

\subsection{Linear convergence of the NIDS with $\hat{g}_i(\cdot)\neq 0$}\label{subsec:linear-rate}

This subsection presents a linear convergence rate of the NIDS via the CD-DYS.
We first present a new result of linear convergence for the general DYS for Problem \eqref{problem_3O} (not limited to the CD-DYS) 
when $f$ is strongly convex and
$g$ is the indicator function of a linear image space.
As indicator functions satisfy neither smoothness nor strong convexity,
our result cannot be derived from the prior results of linear convergence as \cite{davis2017three,condat2022randprox,lee2022convergence,yi2022convergence}.
The proof is presented in Section \ref{sec:proof-rate}.
\begin{theorem}\label{lemma:convergence}
Consider the variable metric DYS in \eqref{alg:DYS_variable_metric_pre} for $k=1,2,\ldots$ for Problem \eqref{problem_3O}.
Let $y^*$ and $z^*$ be the optimal values of $y^{k+1/2}$, $y^k$, and $z^k$.
Suppose that $M^{-1}\nabla f(y)$ is $L$-Lipschitz continuous,
$f$ is $\mu$-strongly convex, $g(y)= \delta_{\mathrm{Im}(U)}(y)$
with a column full-rank matrix $U$, and $h$ is proper, closed, and convex.
Set a stepsize $\alpha \in(0,2\varepsilon/L)$, where $\varepsilon \in(0,1)$.
Pick any start point $y^{1/2} = y^{\mathrm{init}}$ and set $z^1=y^{1/2}-\alpha M^{-1}\nabla f(y^{1/2})$.
Then it holds 
that
\begin{align*}
&\|z^{k+1}-z^*\|^2
+ \nu\|\zeta(y^{k+1/2}) -\zeta(y^{*})\|^2\\
\leq & (1-C)( \|z^{k}-z^*\|^2
+ \nu\|\zeta(y^{k-1/2}) -\zeta(y^{*})\|^2),
\end{align*}
where $\nu\in (0, \frac{\beta}{2}(\alpha-\frac{\alpha^2 L}{2\varepsilon}))$
with $\beta=\min\{\frac{1}{\alpha L},\mu\}$,
$\zeta(y):=y-\alpha M^{-1}\nabla_y f(y)$, and
    $C = \min\{\frac{\kappa}{48},\frac{\kappa}{12\alpha},\frac{\nu}{\nu+9}\}$
with $\kappa:=\beta(\alpha-\frac{\alpha^2 L}{2\varepsilon})-2\nu>0$.
\end{theorem}

Since $\CD_l = \mathrm{Im}(\mathbf{1}_{|\CCl|})$ for $\CD_l$ in \eqref{D_l},
Theorem \ref{lemma:linear-rate} provides the following linear rate for the NIDS with $\BPhi$.
Although \cite{alghunaim2020decentralized,xu2021distributed} have addressed this case, our result below admits bigger stepsizes due to the arbitrariness of $\varepsilon\in(0,1)$.
\begin{theorem}\label{lemma:linear-rate}
Consider the same assumptions as Proposition \ref{prop:NIDS}.
Further, assume that $\CG$ is connected and that for each $i\in\CN$, $\hat{f}_i(\cdot)$ is $\hat{\mu}_i$-strongly convex.
Set a stepsize $\alpha \in(0,2\varepsilon/\max_{i\in\CN}\hat{L}_i)$, where $\varepsilon \in(0,1)$. 
Pick any start point $\Bx^{0}$ and set
$\Bw^1=\Bx^{0}-\alpha \hat{f}(\Bx^{0})$.
Then, $\|\Bx^k-\Bx^*\| = O((1-C)^{k/2})$ holds, where $C$ is given as Theorem \ref{lemma:convergence} with $L=\max_{i\in\CN}\hat{L}_i$ and
$\mu=\min_{i\in\CN}\hat{\mu}_i/\max_{i\in\CN}|\QiGactive|$.
\end{theorem}
\begin{proof}
This theorem follows from Proposition \ref{prop:NIDS} and Theorem \ref{lemma:convergence}.
Note that while the $\mu_i$-strong convexity of $\hat{f}_i$ for $i=1,\ldots,n$ guarantees the convexity of
$\hat{f}((\BD^\top\BD)^{-1}\BD^\top\By) - (\min_i\hat{\mu}_i )\|(\BD^\top\BD)^{-1}\BD^\top\By\|^2$\footnote{ $\hat{f}((\BD^\top\BD)^{-1}\BD^\top\By)$ is strongly convex over $\mathrm{Im}(\BD)$; recall our formulation in \eqref{problem_y_aw}.}, we can treat $\sum_{i=1}^n \hat{f}_i(E_i(\BD^\top\BD)^{-1}\BD^\top\By)$ just as a $ \min_{i\in\CN}\hat{\mu}_i/\max_{i\in\CN}|\QiGactive|$-strongly convex function and
directly apply the same arguments as Theorem \ref{lemma:convergence} because both $\By^{k+1/2}$ in Algorithm \eqref{alg:DYS_distributed_aw_vm} and $\By^*$ always belong to $\mathrm{Im}(\BD)$.
The linear convergence of $\{\Bx^k\}$ follows from the first line of Algorithm \eqref{alg:DYS_distributed_aw_vm}.
\end{proof}

\section{Numerical Experiments}\label{sec:numerical_experiments}
This section presents two numerical examples: resource allocation and consensus optimization. 
\paragraph{Clique-wise resource allocation}
First, we consider a resource allocation problem for the network of $n=20$ agents in \cite[Fig. 1]{watanabe2022distributed} with four communities modeled by cliques and, suppose that,
for the local consumption, a resource constraint is imposed on each community.
The clique parameters are given as
$\QGactive=\{1,2,3,4\}$ and $\CC_1=\{1,2,\ldots,6\},\,\CC_2=\{5,6,\ldots,9\},\,\CC_3=\{8,9,\ldots,12\},\,\CC_4=\{9,10,13,14,\ldots,20\}$.
Let $x_i\geq 0$ for $i\in\CN$ be the amount of resources of agent $i$.
For $l\in\QGactive$,
we set 
the clique-wise cost function 
that can account for the desired resource amount with a weight for each community
as
$f_l(x_\CCl) = \frac{a_l}{2}\| \frac{1}{|\CCl|}\sum_{j\in\CCl} x_j - b_l  \|^2$, with weights $a_l=l,\,l\in\{1,\ldots,4\}$ and desired resource amounts $b_l$ (generated by the uniform distribution with $[0,50]$), 
and define the clique-wise resource constraint 
 as $g_l(x_\CCl) = \delta_{\CD_l}(x_\CCl)$ with $\CD_l=\{x_\CCl:\sum_{j\in\CCl}x_j \leq N_l\}$
for $(N_1,\ldots,N_4)=(5,10,5,15)$.
For $i\in\CN$,
we consider the agent's utility
$\hat{f}_i(x_i)=\frac{\hat{a}_i}{2}\|x_i-\hat{b}_i\|^2$ with $\hat{a}_i=1$, $\hat{b}_i$ generated by the uniform distribution with $[0,10]$, and nonnegativity of consumption $\hat{g}_i(x_i) = \delta_{\BR_+\cup\{0\}}(x_i)$.

We here compare the proposed methods in \eqref{alg:DYS} (or Algorithm \ref{alg:DYS_distributed_aw}) with the parameter $\alpha=1/(\max_{i\in\CN,l\in\QGactive}\hat{a}_i/\min_{i\in\QGactive}|\QiGactive|  +\max_{l\in\QGactive}a_l)$
with Liang et al. \cite{liang2019distributed} for two different stepsizes ($\tau=0.1,\,0.2$).
The method in \cite{liang2019distributed} is a distributed algorithm for globally coupled constraints using the gradient of the cost function.
Notice that the dual decomposition technique cannot be directly used here since we have to minimize
$f_l,\,l\in\QGactive$.

\begin{figure}[t]
    \centering
    \includegraphics[width=0.73\columnwidth]{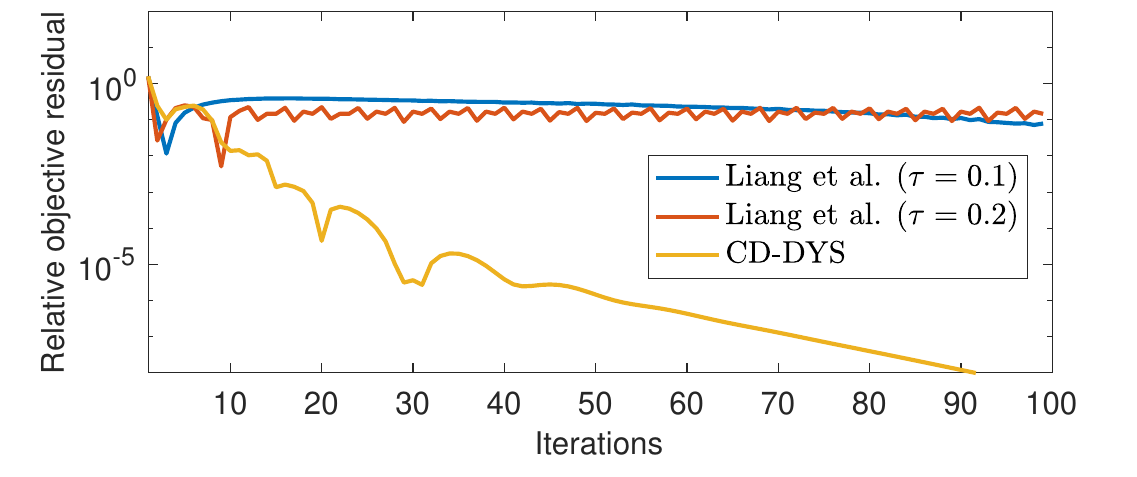}
    \caption{Plots of the relative objective residual under the 
    Liang et al. \cite{liang2019distributed} with $\tau=0.1,\,0.2$ and the CD-DYS in Algorithm \ref{alg:DYS_distributed_aw} (or \eqref{alg:DYS}). Here, $\tau$ represents the step-size.}
    \label{fig:sim_ra}
\end{figure}

The simulation result is plotted in Fig. \ref{fig:sim_ra}.
The proposed CD-DYS converges much faster than Liang et al. \cite{liang2019distributed} whereas that with $\tau=0.2$ fails to give a descent direction.
This difference is rooted in the fact that the CD-DYS exploits the community structure of the problem and admits larger stepsizes.
Note that to get the largest stepsize for Liang et al.,
one has to know an upper bound of the norm of~dual~variables.

\paragraph{Consensus optimization via NIDS}
We next consider the $\ell_1$ norm regularized consensus optimization problem \eqref{problem_clique-wise_constraints} for $n=50$ agents over an undirected graph $\CG$.
Here, we set $\hat{f}_i(x_i) = \frac{1}{2} \|\Psi_i x_i - b_i\|^2$ and $\hat{g}_i(x_i) = \lambda_i \|x_i\|_1$ for $i\in\CN$.
Here, $\Psi_i = I_{10} + 0.05 \Omega_i \in \BR^{10\times 10}$, $b_i\in\BR^{10},\,i\in\CN$, and $\lambda_1=\cdots=\lambda_n=\lambda=0.001$.
Each entry of $\Omega_i$ and $b_i$ is generated by the standard normal distribution.

We here conduct simulations of the NIDS in \eqref{NIDS_prox}
for $\widetilde{\BW} = \BPhi$ with $\QGactive=\maxQG$ and $\QGactive=\QGactive^\mathrm{edge}$, $\widetilde{\BW}_\BL$ with $\epsilon=0.99/\max_{i\in\CN}(|\CN_i|-1)$, $\widetilde{\BW}_\mathrm{mh}$, and $\widetilde{\BW}_c:=I-\alpha c(I-\BW_\BL)$, where $c=1/(1-\lambda_{\min}(\BW_\BL)\alpha)$.
The last is introduced in \cite{li2019decentralized} as a less conservative choice using global information $\lambda_\mathrm{min}(\BW_\BL)$.
The stepsize is assigned as $\alpha=1/\hat{L}$ with $\hat{L}=\max_{i}\{\lambda_\mathrm{max}(|\QiGactive|(\Psi_i^\top\Psi_i))\}.$

\begin{figure}[t]
    \centering
    \includegraphics[width=0.7\columnwidth]{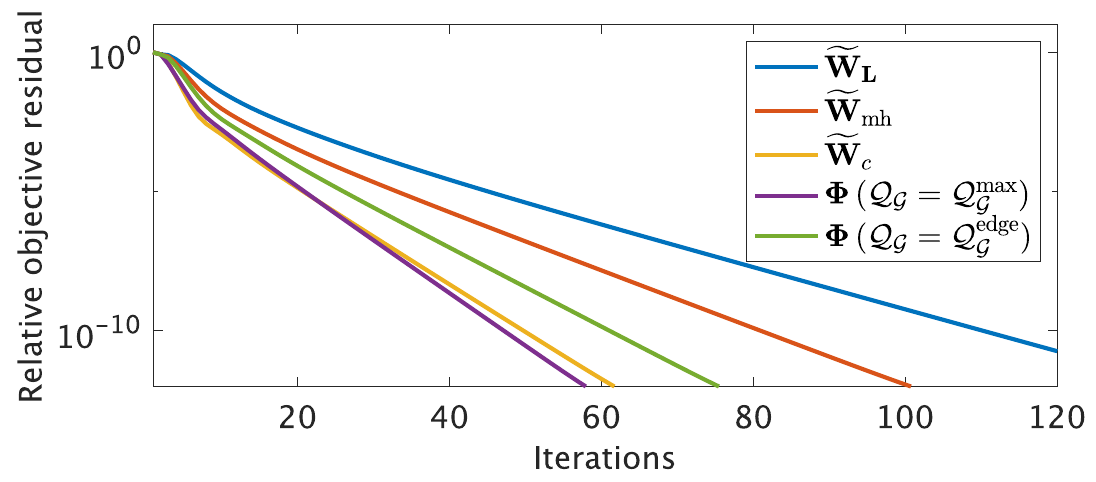}
    \caption{Plots of the relative objective residual under the NIDS in \eqref{NIDS_prox} for the five different choices of $\widetilde{\BW}$.}
    \label{fig:sim_consensus}
\end{figure}

The simulation result is reported in Fig. \ref{fig:sim_consensus}, which plots the relative objective residual $|F(\Bx^k)-F(\Bx^*)|/F(\Bx^*)$ where $F(\Bx):=\hat{f}(\Bx)+\lambda\|\Bx\|_1$.
It can be observed that the case of $\BPhi$ with $\QGactive=\maxQG$ exhibits the fastest convergence in almost 60 iterations with high accuracy, outperforming $\widetilde{\BW}_c$ without using global information.
While the case of $\BPhi$ with $\QGactive=\QGactive^\mathrm{edge}$ is slower than $\widetilde{\BW}_c$, this is still superior to $\widetilde{\BW}_\mathrm{mh}$ and $\widetilde{\BW}_\BL$, as implied in Proposition \ref{prop:Phi_eigval}.


\section{Proof of Theorem \ref{lemma:linear-rate}}\label{sec:proof-rate}
Our proof is based on the following trick (See Theorem D.6 in \cite{davis2017three}):
If $a_0, \cdots, a_p, b_0, \cdots, b_p, c_0, \cdots, c_p \in (0,\infty) $ for some $p > 0$, and
\begin{equation}\label{rate_key_inequality}
\text{{${\!\!\!\!\!\|w^{k+1}-w^*\|^2+\sum_{i=0}^p a_i c_i \leq\|w^k-w^*\|^2 \leq \sum_{i=0}^p a_i b_i}$}},
\end{equation}
then $\sum_{i=0}^n a_i b_i \leq \max_{i} (b_i / c_i)\sum_{i=0}^n a_i c_i$, so
\begin{align*}
\|w^{k+1}-w^*\|^2+\min_{i} (c_i / b_i) \|w^k-w^*\|^2 
\leq\|w^{k+1}-w^*\|^2+\sum_{i=0}^p a_i c_i \leq\|w^{k}-w^*\|^2 \leq\|w^{k}-w^*\|^2.
\end{align*}
Thus,
\begin{align}\label{rate_trick}
    \!\!\!\!\!\!
    \|w^{k+1}-w^*\|^2 \leq(1-\min_{i} c_i / b_i )\|w^{k}-w^*\|^2,
\end{align}
which provides a linear convergence rate.
In the following, we derive an inequality of the form in \eqref{rate_key_inequality}
with
$w^k = [(z^k )^\top,
\nu^{1/2}(\zeta(y^{k-1/2}))^\top]^\top
$, $w^* = [z^{*\top},
\nu^{1/2}(\zeta(y^{*}))^\top]^\top
$, and
some constants $a_0, \cdots, a_p, b_0, \cdots, b_p, c_0, \cdots, c_p>0$.

We first prepare a key inclusion for establishing the desired rate.
We suppose that $z^k$, $y^k$, and $y^{k+1/2}$ are not optimal without loss of generality.
For $g(y)= \delta_{\mathrm{Im}(U)}(y)$, we have
$\prox_{\alpha g}^M(y)
=U(U^\top M U)^{-1}U^\top M y$, which leads to
\begin{align*}
    U^\top M z^{k+1}
    &= U^\top M z^{k} -IMy^{k+1/2}
    +U^\top M(2y^{k+1/2} - z^k -\alpha M^{-1} \nabla_y f(y^{k+1/2}) )\\
    &= U^\top M(y^{k+1/2} - \alpha M^{-1}\nabla_y f(y^{k+1/2}))
\end{align*}
since $U^\top M \prox_{\alpha g}^M(y) = U^\top My $.
Note that this holds for all $k\geq 1$ thanks to the initialization.
Then we have
\begin{align*}
y^{k+1} = (U^\top MU)^{-1} U^\top M (2y^{k+1/2} - y^{k-1/2}
-\alpha M^{-1} \nabla_y f(y^{k+1/2})
+\alpha M^{-1} \nabla_y f(y^{k-1/2})),
\end{align*}
and thus we get $\xi^{k+1}:=2y^{k+1/2} - y^{k-1/2}
-\alpha M^{-1} \nabla_y f(y^{k+1/2})
+\alpha M^{-1} \nabla_y f(y^{k-1/2})
\in (I+\alpha M^{-1} \pd g) (y^{k+1})$
by \cite[Prop. 16.44]{bauschke2017correction}.
This inclusion allows us to remove the assumption of smoothness or strong convexity for $g$ or $h$ in \cite{davis2017three,lee2022convergence,yi2022convergence,condat2022randprox} because one can evaluate $y^{k+1}$ only with $f$,
without $z^k$.

We derive the lower side of the inequality in \eqref{rate_key_inequality} as follows.
By the smoothness and strong convexity of $f$,
for $\|z^{k}-z^*\|^2$, \cite[Proposition D.4]{davis2017three} provides
    \begin{align*}
\|z^{k}-z^*\|^2
\geq&(1-\varepsilon)\|z^{k+1}-z^*\|^2
+ 2\alpha \max\{\mu \|y^{k+1/2} -y^*\|^2, 
\frac{1}{L}\|M^{-1}\nabla_y f(y^{k+1/2})-M^{-1}\nabla_y f(y^{*})\|^2\}\\
&-\frac{\alpha^2}{\varepsilon}\|M^{-1}\nabla_y f(y^{k+1/2})-M^{-1}\nabla_y f(y^{*})\|^2.
\end{align*}
Then, using $\theta=1 /(2-\varepsilon)$, $\|w^k-w^*\|^2$ is lower bounded as \begin{align}\label{convergence_rate_lower}
&\|w^k-w^*\|^2
\geq \|z^{k+1} - z^*\|^2 +\nu \|\zeta(y^{k-1/2})-\zeta(y^*) \|^2\nonumber
\\
&+(\frac{1}{\theta}-1)\|z^{k+1}-z^k\|^2 
+2\alpha \max\{\mu \|y^{k+1/2}-y^*\|^2,
\frac{1}{L}\|M^{-1}\nabla_y f(y^{k+1/2})-M^{-1}\nabla_y f(y^{*})\|^2\} \nonumber\\
&- \frac{\alpha^2 L}{\varepsilon}\times \frac{1}{L}\|M^{-1}\nabla_y f(y^{k+1/2})-M^{-1}\nabla_y f(y^{*})\|^2\nonumber\\
\geq&
\|w^{k+1}-w^*\|^2
+c_0\|z^{k+1}-z^k\|^2
+c_1 \|y^{k+1/2}-y^*\|^2
+c_2 \|M^{-1}\nabla_y f(y^{k+1/2})-M^{-1}\nabla_y f(y^{*})\|^2 
\nonumber\\
&+c_3 \|\zeta(y^{k-1/2})-\zeta(y^*) \|^2,
\end{align}
where
$c_0=1/\theta -1$, $c_1=\beta(\alpha-\frac{\alpha^2L}{2\varepsilon})-2\nu$, $c_2 = \alpha c_1$, and $c_3=\nu$.
Here, the term $\|w^{k+1}-w^*\|^2$ in
\eqref{convergence_rate_lower} comes from the definition of $w^k$ and
the relationship $\|w\|^2+\|w'\|^2 \geq \|w+w'\|^2/2$ obtained by Jensen's inequality for the terms of $\|y^{k+1/2}-y^*\|^2$ and $\|M^{-1}\nabla_y f(y^{k+1/2})-M^{-1}\nabla_y f(y^{*})\|^2$.

Next, utilizing a similar calculation to the proof of \cite[Proposition D.4]{davis2017three} (see the proof of the second upper bound)
and the inclusion $\xi^{k+1}\in(I+\alpha M^{-1}\pd g)(y^{k+1})$,
we derive the upper-bound of $\|z^{k}-z^*\|^2 + \nu \|\zeta(y^{k-1/2})-\zeta(y^*) \|^2$ as follows.
Here, we set $\epsilon = c_0/C$, and
let $u_A^{k+1} \in M^{-1}\pd g(y^{k+1})$ and $u_A^{*} \in M^{-1}\pd g(y^{*})$ satisfy $y^{k+1}+\alpha u_A^{k*1} = \xi^{k+1}$ and
$y^{*}+\alpha u_A^{*} = \xi^{*}$, respectively:
\begin{align}\label{convergence_rate_upper}
&\|w^k-w^*\|^2
\leq 
\nu \|\zeta(y^{k-1/2})-\zeta(y^*) \|^2 
+\|y^{k+1} - \alpha (u_A^{k+1}
+ M^{-1}\nabla_y f(y^{k+1/2})) \nonumber\\
&+ 2(y^{k+1/2}-y^{k+1})
-(y^*-\alpha u_A^*-\alpha \nabla_y f(y^*))
\|^2\nonumber\\
\leq &
\nu \|\zeta(y^{k-1/2})-\zeta(y^*) \|^2 +\|-(\xi^{k+1}-\xi^*)
- \alpha(M^{-1}\nabla_y f(y^{k+1/2}) - M^{-1}\nabla_y f(y^{*})) \nonumber\\
&+2(y^{k+1/2}-y^{*})
\|^2 + \epsilon \|z^{k+1}- z^{k}\|^2 \nonumber\\
\leq &
\nu \|\zeta(y^{k-1/2})-\zeta(y^*) \|^2  + \epsilon \|z^{k+1}- z^{k}\|^2
+3\|\xi^{k+1}-\xi^*\|^2 + 12\|(y^{k+1/2}-y^{*})\|^2\nonumber\\
&+ 3\alpha^2\|M^{-1}\nabla_y f(y^{k+1/2}) - M^{-1}\nabla_y f(y^{*})\|^2 \nonumber\\
\leq &
(\nu+9) \|\zeta(y^{k-1/2})-\zeta(y^*) \|^2  + \epsilon \|z^{k+1}- z^{k}\|^2\nonumber\\
&+ 48\|(y^{k+1/2}-y^{*})\|^2 
+ 12\alpha^2\|M^{-1}\nabla_y f(y^{k+1/2}) - M^{-1}\nabla_y f(y^{*})\|^2,
\end{align}
where the last line follows from
$\|\xi^{k+1}-\xi^*\|^2 \leq 3\|\zeta(y^{k-1/2})-\zeta(y^*) \|^2+12\|(y^{k+1/2}-y^{*})\|^2+3\alpha^2\|M^{-1}\nabla_y f(y^{k+1/2}) - M^{-1}\nabla_y f(y^{*})\|^2$.

Therefore, for $a_0=\|z^{k+1}- z^{k}\|^2$, $a_1=\|(y^{k+1/2}-y^{*})\|^2$, $a_2=\|M^{-1}\nabla_y f(y^{k+1/2}) - M^{-1}\nabla_y f(y^{*})\|^2$, $a_3=\|\zeta(y^{k-1/2})-\zeta(y^*) \|^2$,
$b_0=\epsilon$, $b_1=48$, $b_2=12\alpha^2$, and $b_3=\nu+9$,
the linear rate follows from \eqref{rate_trick}, \eqref{convergence_rate_lower}, and \eqref{convergence_rate_upper}.

\section{Conclusion}\label{sec:conclusion}

This note addressed distributed optimization of clique-wise coupled problems via operator splitting.
First, we defined the CD matrix and a new mixing matrix and analyzed its properties.
Then, using the CD matrix, we presented the CD-DYS algorithm via the Davis-Yin splitting (DYS).
Subsequently, its connection to consensus optimization methods as NIDS was also analyzed.
Moreover, we presented a new linear convergence rate not only for the NIDS with non-smooth terms but also for the general DYS with a projection onto a subspace.
Finally, we demonstrated the effectiveness via numerical examples.

\bibliographystyle{unsrt}
\bibliography{bibliography.bib}

\newpage
\appendix
\section{Application examples}\label{A:examples}

In addition to the examples in Section \ref{sec:numerical_experiments},
we here present various application examples of Problem \ref{problem}.
\paragraph{Formation control}

A formation control problem aims to steer the positions of robots to a desired configuration and has been actively investigated for the past two decades.
For a multi-agent system over undirected graph $\CG=(\CN,\CE)$, the most basic formulation of this problem is 
\begin{align}\label{problem:formation}
\begin{array}{cl}
\underset{
\substack{x_i\in \BR^{d_i},\,i\in\CN}
}
{\mbox{minimize}}&
\displaystyle  \sum_{\{i,j\}\in\CE} \|x_i-x_j - d_{ij}\|^2,
\end{array}
\end{align}
where $x_i$ is the position of agent $i$, and $d_{ij}$ is the desired relative position from $x_j$ to $x_i$.
By assigning $\QGactive=\QGactive^\mathrm{edge}$, one can obtain the desired configuration via the proposed CD-DYS.
Note that one can also deal with various constraints in the clique-wise coupled framework, e.g., an agent-wise constraint $x_i\in\Omega_i$ and a pairwise distance constraint $ \underline{\delta}_{ij}  \leq\|x_i-x_j\|\leq \overline{\delta}_{ij}$.
In addition, the proposed framework also allows us to achieve the desired formation in a distributed manner even for linear multi-agent systems, as shown in \cite{latafat2019new}, and in the case where each agent has no access to the global coordinate and can only use information via relative measurements, as shown in \cite{sakurama2020unified,sakurama2022generalized}
\paragraph{Network Lasso}

The network lasso is an optimization-based machine-learning technique accounting for network structures.
For a multi-agent system over graph $\CG=(\CN,\CE)$, a network lasso problem \cite{hallac2015network} is given as follows:
\begin{align}\label{problem:network_lasso}
\begin{array}{cl}
\underset{
\substack{x_i\in \BR^{m},\,i\in\CN}
}
{\mbox{minimize}}&
\displaystyle  \sum_{i=1}\hat{f}_i(x_i) + \lambda  \sum_{\{i,j\}\in\CE} w_{ij}\|x_i-x_j\|,
\end{array}
\end{align}
where $\lambda>0$ and $w_{ij}>0$ for $\{i,j\}\in\CE$.
This problem can be seen as a special case of Problem \eqref{problem}.
Owing to the second term in \eqref{problem:network_lasso}, neighboring nodes are more likely to form a cluster, i.e., to take close values.
Applications of the Network Lasso include the estimation of home prices \cite{hallac2015network}, where there is a spatial interdependence among houses' prices.

\paragraph{Sparse semidefinite programming}

\begin{figure}
    \centering
    \includegraphics[width = 0.6\columnwidth]{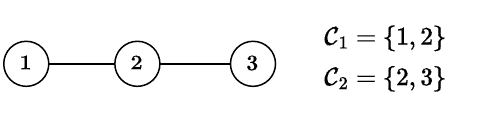}
    \caption{Example of a system with three nodes in Example \ref{ex:CDmatrix}.}
    \label{fig:CDmatrix_ex}
\end{figure}

Semidefinite programming via chordal graph decomposition has been actively studied not only in optimization \cite{vandenberghe2015chordal,fukuda2001exploiting} but also in control \cite{zheng2021chordal} as an efficient and scalable computation scheme exploiting the sparsity of matrices that naturally arises from underlying networked structures of problems.
This type of problem can also be solved in a distributed manner based on the framework of clique-wise coupling.

Consider a multi-agent system over $\CG=(\CN,\CE)$ and the following standard semidefinite programming
\begin{align}\label{problem:SDP}
\begin{array}{cl}
\underset{
\substack{y_i\in\BR,\,i\in\CN,\,Z}
}
{\mbox{minimize}}&
\displaystyle  \sum_{i=1}^n b_i^\top y_i  + \delta_{\mathbb{S}^n_+}(Z)
\\
{\mbox{subject to}}&  \displaystyle Z + \sum_{j=1}^p A_j y_j =C,
\end{array}
\end{align}
where we consider that agent $i$ possesses $y_i$ and $i$th column of $Z$.
Here, $\mathbb{S}_+^n$ represents the set of $n\times n$ positive semidefinite matrices.
This problem cannot be solved in a distributed manner by standard algorithms due to the undecomposable constraint $Z\in\mathbb{S}^n_+$.
Nevertheless, if $Z,\,A_1,\ldots,A_p,\,C$ have the sparsity with respect to $\CG=(\CN,\CE)$ and graph $\CG$ is \textit{chordal},
\cite{vandenberghe2015chordal,zheng2021chordal} show that this problem can be equivalently transformed into the following decomposed form with smaller positive semidefinite constraints: 
\begin{align}\label{problem:SDP_decomp}
\begin{array}{cl}
\underset{
\substack{y_i\in\BR,\,i\in\CN,\,Z_l,\,l\in\maxQG} 
}
{\mbox{minimize}}&
\displaystyle  \sum_{i=1}^n b_i^\top y_i  + \sum_{l\in\maxQG} \delta_{\mathbb{S}_+^{|\CCl|}} (Z_l)\\
{\mbox{subject to}}&  \displaystyle  \sum_{j=1}^p A_j y_j + \sum_{l\in\maxQG}D_l^\top Z_l D_l=C.
\end{array}
\end{align}
Moreover, when $\sum_{i=1}b_i^\top y_i $ in \eqref{problem:SDP_decomp} can be rewritten as 
\begin{equation}\label{MYYN}
    \sum_{j=1}^p A_j y_j = M Y +Y N
\end{equation}
with $Y= \diag(y_1,\ldots,y_n)$ and some matrices $M,N$ with the sparsity with respect to $\CG=(\CN,\CE)$, we can reformulate Problem \eqref{problem:SDP_decomp} into a clique-wise coupled problem in \eqref{problem} by introducing auxiliary variables.
For example, for the system with $n=3$ in Fig. \ref{fig:CDmatrix_ex}, which is a chordal graph with maximal cliques $\CC_1=\{1,2\}$ and $\CC_2=\{2,3\}$, Problem \eqref{problem:SDP_decomp} with \eqref{MYYN} reduces to 
\begin{align}\label{problem:SDP_decomp_example}
\begin{array}{cl}
&\underset{
\substack{y_1,y_2,y_3\in\BR,\,Z_1,Z_2\in\mathbb{S}^2} 
}
{\mbox{minimize}}\quad
\displaystyle  \sum_{i=1}^3 b_i^\top y_i  + \delta_{\mathbb{S}_+^{2}} (Z_1)+\delta_{\mathbb{S}_+^{2}} (Z_2)\\
&{\mbox{subject to}}\\
&\begin{bmatrix}
m_{11}y_1 +n_{11}y_1 & m_{12}y_1+n_{12}y_2 &0 \\
m_{21}y_2 + n_{21}y_1& m_{22}y_2+n_{22}y_2 &m_{23} y_2+n_{23}y_3\\
0 & m_{32}y_3+n_{32}y_2 &m_{33} y_3 + n_{33}y_3
\end{bmatrix}     
+ \begin{bmatrix}
\text{\large{$Z_1$}}
& \begin{matrix}
    0 \\0
\end{matrix}\\
\begin{matrix}
    0 &0
\end{matrix}&0
\end{bmatrix} 
+ \begin{bmatrix}
0&\begin{matrix}
    0 &0
\end{matrix} \\
 \begin{matrix}
    0 \\0
\end{matrix}& \text{\large{$Z_2$}}
\end{bmatrix} 
=\begin{bmatrix}
c_{11} & c_{12} &0 \\
c_{21}& c_{22} &c_{23}\\
0 & c_{32} &c_{33}
\end{bmatrix},
\end{array}
\end{align}
where $m_{ij}$, $n_{ij}$, $c_{ij}$ are the $i,j$ entries of $M$, $N$, and $C$, respectively.
Hence,
by decomposing the constraint in \eqref{problem:SDP_decomp_example} into clique-wise coupled constraints by using the auxiliary variables $ \hat{z}_{2,11}$ and $ \hat{z}_{1,22}$
as
\begin{align}
\label{sdp_decomposed_constraints_1}
&\begin{bmatrix}
m_{11}y_1 +n_{11}y_1 & m_{12}y_1+n_{12}y_2 \\
m_{21}y_2 + n_{21}y_1& m_{22}y_2+n_{22}y_2 
\end{bmatrix} 
+ Z_1 
+ \begin{bmatrix}
0 & 0\\
0 & \hat{z}_{2,11}
\end{bmatrix}
= \begin{bmatrix}
c_{11} & c_{12}\\
c_{21}& c_{22} \\
\end{bmatrix}   \\
\label{sdp_decomposed_constraints_2}
&\begin{bmatrix}
m_{22}y_2 +n_{22}y_2 & m_{23}y_2+n_{23}y_3 \\
m_{32}y_3 + n_{32}y_2& m_{33}y_3+n_{33}y_3 
\end{bmatrix} 
+ Z_2
+ \begin{bmatrix}
\hat{z}_{1,22} & 0\\
0 &0
\end{bmatrix}
= \begin{bmatrix}
c_{22} & c_{32}\\
c_{32}& c_{33} \\
\end{bmatrix}    \\
\label{sdp_decomposed_constraints_3}
&z_{1,22} = \hat{z}_{1,22},\quad z_{2,11} = \hat{z}_{2,11},
\end{align}
where $z_{l,ij}$ represents the $i,j$ entry of $Z_l$,
we can obtain an equivalent clique-wise coupled problem in the following with 
$x_1 = [y_1; \mathrm{vec}([Z_1]_1)]$, $x_2 = [y_2;\mathrm{vec}([Z_1]_2);\mathrm{vec}([Z_2]_1);\hat{z}_{1,22};\hat{z}_{2,11}]$, $x_3 = [y3; \mathrm{vec}([Z_2]_2)]$, where $\mathrm{vec}([Z_l]_i)$ represents $i$th column of the matrix $Z_l$:
\begin{align}\label{problem:SDP_decomp_example_2}
\begin{array}{cll}
&\underset{
\substack{y_1,y_2,y_3\in\BR\\Z_1,Z_2\in\mathbb{S}^2} 
}
{\mbox{minimize}}&
\displaystyle  \sum_{i=1}^3 b_i^\top y_i  + \delta_{\mathbb{S}_+^{2}} (Z_1)+\delta_{\mathbb{S}_+^{2}} (Z_2)\\
&{\mbox{subject to}} &\text{\eqref{sdp_decomposed_constraints_1}--\eqref{sdp_decomposed_constraints_3}}
\end{array}
\end{align}
Semidefinite programming in the form of \eqref{problem:SDP_decomp} with \eqref{MYYN} arises in practical problems, e.g., distributed design of decentralized controllers for linear networked systems \cite{zheng2019distributed} and sensor network localization \cite{anjos2011handbook}.
Note that one can extend the discussion above to higher dimensional vectors $y_i$ and block-partitioned matrices $Z$, as shown in \cite{zheng2021chordal}.

\paragraph{Approximating trace norm minimization problems}

Trace norm minimization is a powerful technique in machine learning and computer vision that can obtain a low-rank matrix $\hat{L}$ representing the underlying structure of the data.
Its applications include the Robust PCA (RPCA) and multi-task learning problems.

For example, we can relax an RPCA problem to a clique-wise coupled problem as follows.
Consider a data matrix $Y\in\BR^{d\times n}$.
Then, a standard form of RPCA is formulated as follows:
\begin{align}
\label{prob:rpca}
\underset{\hat{S},\hat{L}}{\mbox{minimize}}\;\;
\displaystyle \|\hat{S}\|_1 + \theta \|\hat{L}\|_*  \;\;
\displaystyle\mbox{subject to}\;\; \hat{S}+\hat{L} = Y.
\end{align}
By solving this problem, we can decompose a data matrix $Y$ into two components: a low-rank matrix $\hat{L}$ representing the underlying structure of the data and a sparse matrix $\hat{S}$ capturing the outliers or noise.
Consider that for a multi-agent system with $n$ agents and
$Y=[Y_1,\ldots,Y_n]$, agent $i$ possesses the matrix $Y_i$.
Then the robust PCA problem with the clique-based relaxation is formulated as follows:
\begin{align}\label{prob:rpca_clique}
\begin{array}{cl}
\underset{\substack{\hat{S}_i,\, i\in\CN \\ 
\hat{L}_l \in \BR^{640\times 40|\CCl|},\,l\in\QGactive }}{\mbox{minimize}}&
\displaystyle \sum_{i\in\CN} \|\hat{S}_1\|_1 + \sum_{l\in\CCl} \theta_l \|\hat{L}_l\|_*  \quad \\
\displaystyle\mbox{subject to} & 
\hat{S}_\CCl + \hat{L}_l = Y_\CCl \quad \forall l\in\QGactive,
\end{array}
\end{align}
where $\hat{S}_\CCl= [\hat{S}_{j_1},\ldots,\hat{S}_{j_{|\CCl|}}]$ and $\hat{Y}_\CCl= [\hat{Y}_{j_1},\ldots,\hat{Y}_{j_{|\CCl|}}]$ for $\CCl=\{j_1,\ldots,j_{|\CCl|}\}$.
Here, $\hat{S}_i$ and $\hat{L}_l$ correspond to $x_i$ and $y_l$ in Problem \eqref{problem}.
Although Problem \eqref{prob:rpca_clique} involves relaxation, one can still realize a low-rank matrix by solving it.

\section{Proof of Lemma \ref{thorem:CDmatrix}}
(a) We prove the statement by contradiction. Assume that the CD matrix $\BD$ is not column full rank. Then, there exists a vector $\Bv=[v_1^\top,\ldots,v_n^\top]^\top\neq 0$ with $v_i\in\BR^{d_i}$ such that $\BD \Bv =0$. This yields $D_l \Bv =0$ for $\Bv$ and all $l\in\QGactive$. Hence, we obtain $E_i \Bv = v_i = 0$ for all $i\in\CN$ from Assumption \ref{assumption:sumCl=N}.
This contradicts the assumption.

(b) For $\BD$, we have $\BD^\top \BD = \sum_{l\in\QGactive} D_l^\top D_l = \sum_{l\in\QGactive}\sum_{j\in\CCl} E_j^\top E_j = \sum_{i=1}^n \sum_{l\in \QiGactive} E_i^\top E_i = \sum_{i=1}^n |\QiGactive| E_i^\top E_i  $ from Definition \ref{def:CDmatrix}.
Here, $E_i^\top E_i = \bdiag(O_{d_1\times d_1},\ldots, I_{d_i},\ldots,O_{d_n\times d_n})$ holds.
Therefore, we obtain $\BD^\top \BD = \bdiag (|\CQ_\CG^1| I_{d_1},\ldots,|\CQ_\CG^n| I_{d_n})$.
$\BD^\top\BD \succ O$ follows from Assumption \ref{assumption:sumCl=N}.

(c) It holds that $\BD^\top \By = \sum_{l\in\QGactive} D_l^\top y_l = \sum_{l\in\QGactive} \sum_{j\in\CCl} E_j^\top (E_{l,j} y_l) = \sum_{i=1}^n \sum_{l\in\QiGactive} E_i^\top E_{l,i} y_l = \sum_{i=1}^n E_i^\top (\sum_{l\in\QiGactive} E_{l,i} y_l )$. 
Hence, we obtain \eqref{eq:DTtimesy}.
\section{Proof of Proposition \ref{thorem:CDmatrix_grad_prox}}
(a) For $\Bz\in\mathrm{Im}(\BD)$, there exists some $\Bx\in\BR^d$ such that $\Bz = \BD \Bx$. Then, we obtain
\begin{align*}
    &\prox_{\alpha G}(\By)
    = \BD  \argmin_{\Bx\in\BR^d} ( \frac{1}{2\alpha} \|\By - \BD\Bx\|^2 + \sum_{i=1}^n \hat{g}_i(E_i \Bx)) \nonumber \\
    &= \BD \argmin_{\Bx\in\BR^d} ( \sum_{i=1}^n  (\sum_{l\in\QiGactive} \frac{1}{2\alpha}\|E_{l,i}y_l  - x_i\|^2 +  \hat{g}_i(x_i))) \nonumber \\
    &= \BD   \argmin_{\Bx\in\BR^d} ( \sum_{i=1}^n ( \frac{|\QiGactive|}{2\alpha} \|\sum_{l\in\QiGactive}\frac{1}{|\QiGactive|} E_{l,i}y_l  - x_i\|^2 +  \hat{g}_i(x_i))).
\end{align*}
Therefore, we obtain \eqref{eq:prox_CDmatrix}. Note that the last line can be verified by considering the optimality condition.

(b) This can be proved in the same way as Proposition \ref{thorem:CDmatrix_grad_prox}a with an easy modification from the definition of $\BQ$.

(c) By the chain rule, we have
    $\frac{\pd}{\pd \By}\hat{f}_i(E_i (\BD^\top\BD)^{-1}\BD^\top \By ) = \BD (\BD^\top \BD)^{-1} E_i^\top \nabla_{x_i} \hat{f}_i(E_i (\BD^\top\BD)^{-1}\BD^\top \By )$.
which gives \eqref{eq:grad_CDmatrix}.

\section{Proof of Proposition \ref{prop:CDmatrix_variable_metric}}

(a) For $\BQ$, we obtain $\BQ\BD = [Q_l D_l]_{l\in\QGactive}$.
Then, 
\begin{equation*}
    \BD^\top\BQ\BD = \sum_{l\in\QGactive} D_l^\top Q_l D_l= \sum_{l\in\QGactive}\sum_{j\in\CCl} \frac{1}{|\QjGactive|} E_j^\top E_j.
\end{equation*}
Thus, following the same calculation as the proof of Lemma \ref{thorem:CDmatrix}b gives $\BD^\top\BQ\BD = I_d$.

(b) For any $\By=[y_l]_{l\in\QGactive} \in \BR^{\hat{d}}$, it holds that 
\begin{equation*}
    \BD^\top \BQ \By = \sum_{l\in\QGactive} D_l^\top Q_l y_l=  \sum_{l\in\QGactive}\sum_{j\in\CCl} 
 \frac{1}{|\QjGactive|}E_j^\top E_{l,j} y_l.
\end{equation*}
Hence, reorganizing this and using the proof of Lemma \ref{thorem:CDmatrix}c yield
\begin{equation*}
   \BD^\top \BQ \By= \sum_{i=1}^n \frac{1}{|\QiGactive|}E_i^\top \sum_{l\in\QiGactive} E_{l,i} y_l = \bdiag([\frac{1}{|\QiGactive|}I_{d_i}]_{i\in\CN}) \BD^\top \By. 
\end{equation*}
  Therefore, we obtain $\BD^\top \BQ = (\BD^\top\BD)^{-1}\BD^\top$ from Lemma \ref{thorem:CDmatrix}b.
  The latter equation is also proved in the same way.

(c) From Proposition \ref{prop:CDmatrix_variable_metric}b and Assumption \ref{assumption:sumCl=N}, it holds that $\BD^\top = (\BD^\top\BD)^{-1}\BD^\top \BQ^{-1}$. For the transpose of this matrix, multiplying $\BD^\top \BD$ from the right side gives $ \BQ^{-1}\BD= \BD(\BD^\top\BD)$. 
The latter equation is also proved in the same manner.

\section{Connection of the CD-DYS to other distributed optimization algorithms}\label{appendix:relation}

We here present the comprehensive analysis for the diagram in Fig. \ref{fig:relation} by deriving the CPGD algorithm in \cite{watanabe2022distributed} and Section \ref{A:CPGD}.

\paragraph{Exact diffusion and Diffusion algorithms}

Over undirected graphs, the exact diffusion algorithm is just a special case of the NIDS.
In the case of $\hat{g}_i=0$ for all $i\in\CN$, the NIDS reduces to the Exact diffusion \cite{yuan2018exact,yuan2018exact2}, which is given as follows:
\begin{align}\label{NIDS}
\!\!\!\!
\Bx^{k+1} &= \widetilde{\BW} (2\Bx^{k} - \Bx^{k-1} 
+ \alpha (\nabla_\Bx \hat{f}(\Bx^{k-1}) - \nabla_\Bx \hat{f}(\Bx^{k}) ) ).
\!\!\!\!
\end{align}
This can be rewritten as follows:
\begin{align}\label{Exact_diffusion}
\begin{array}{lll}
     \Bv^{k+1}=& \Bx^k - \alpha \nabla_\Bx \hat{f}(\Bx^{k}) \\
    \Bx^{k+1} =& \widetilde{\BW} (\Bv^{k+1} + \Bx^k - \Bv^{k}).
\end{array}
\end{align}
Those algorithms exactly converge to an optimal solution under mild conditions.
Note that the Exact diffusion is also valid for directed networks and non-doubly stochastic $\BW$. For details, see \cite{yuan2018exact,yuan2018exact2}.

The diffusion algorithm \cite{sayed2014diffusion,chen2012diffusion} is an early distributed optimization algorithm, given as
\begin{align}\label{diffusion}
\begin{array}{lll}
    \Bx^{k+1} =& \widetilde{\BW} (\Bx^k - \alpha \nabla_\Bx \hat{f}(\Bx^{k})).
\end{array}
\end{align}
This algorithm is obtained from the NIDS for $\hat{g}_i=0,\,i\in\CN$ and the Exact diffusion approximating $ \Bx^k - \Bv^{k}\approx 0$ in the second line of \eqref{Exact_diffusion}.
Notice that conditions on ${\BW}$ in \eqref{diffusion} are not equivalent to \eqref{NIDS} and \eqref{Exact_diffusion} (see \cite{sayed2014diffusion,chen2012diffusion,ryu2022large,yuan2018exact,yuan2018exact2}).
Although its convergence is inexact over constant $\alpha$, its simple structure allows us to easily apply it to stochastic and online setups.

\paragraph{CPGD generalizes of the diffusion algorithm}

Invoking the relationship between NIDS/Exact diffusion and diffusion algorithms,
we derive a diffusion-like algorithm from the variable metric CD-DYS in \eqref{alg:DYS_distributed_aw_vm} for 
\begin{align}\label{problem_clique-wise_constraints_noprox}
\begin{array}{cl}
\underset{
\substack{x_i\in \BR^{d_i},\,i\in\CN}
}
{\mbox{minimize}}&
\displaystyle  \sum_{i=1}^n \hat{f}_i (x_i)+\sum_{l\in\QGactive} \delta_{\CD_l}(x_\CCl),
\end{array}
\end{align}
where $\CD_l$ is a closed convex set and not limited to \eqref{D_l}.
The derived algorithm will be formalized as the clique-based projected gradient descent (CPGD) in Appendix \ref{A:CPGD}.

We derive the diffusion-like algorithm as follows.
From $\hat{g}_i=0$, we have $\Bx^k = \Bx^{k-} = (\BD^\top\BD)^{-1}\BD^\top \Bz^k$ and $(\BD^\top\BD)^{-1}\BD^\top \times \By^{k+1/2} = \Bx^{k}$.
Accordingly, the variable metric CD-DYS in \eqref{alg:DYS_distributed_aw_vm} 
reduces to
\begin{align*}
    \begin{array}{lll}
         &\Bx^k = (\BD^\top\BD)^{-1}\BD^\top \Bz^k  \\
         &\By^{k+1} = P_{\Pi_{l\in\QGactive}\CD_l}^\BQ(2\BD\Bx^k - \Bz^k - \alpha\BD \nabla_\Bx \hat{f}(\Bx^k))\\
         &\Bz^{k+1}= \Bz^k + \By^{k+1} - \BD\Bx^k.
    \end{array}
\end{align*} 
 By using $\Bv^{k+1}$ of the form in \eqref{Exact_diffusion}, we get
\begin{align}\label{alg:diffusion_pre}
 &\Bv^{k+1} = \Bx^k - \alpha \nabla_\Bx \hat{f}(\Bx^k) \\
 &\Bx^{k+1} = (\BD^\top\BD)^{-1}\BD^\top P_{\Pi_{l\in\QGactive}\CD_l}^\BQ(\BD\Bv^{k+1} + \BD\Bx^k -\Bz^k )\nonumber
\end{align}
with $\Bz^k$ 
from 
    $\Bx^{k+1}=(\BD^\top\BD)^{-1}\BD^\top\Bz^{k+1} = (\BD^\top\BD)^{-1}\BD^\top (\Bz^k +\By^{k+1})-\Bx^k = (\BD^\top\BD)^{-1}\BD^\top \By^{k+1}$.
In consensus optimization, it can be observed from the previous subsection that $P_{\Pi_{l\in\QGactive}\CD_l}^\BQ(\cdot)$ boils down to a linear map and
$\Bz^k$ satisfies $P_{\Pi_{l\in\QGactive}\CD_l}^\BQ(\Bz^k) = P_{\Pi_{l\in\QGactive}\CD_l}^\BQ(\BD\Bv^k)$
because we have
\begin{align*}
P_{\CD_l}^{Q_l}(z_l^{k+1}) 
=  P_{\CD_l}^{Q_l} (x^k_\CCl - \alpha D_l \nabla_\Bx \hat{f}(\Bx^k)) =P_{\CD}^{Q_l}(D_l\Bv^{k})
\end{align*}
for $\CD_l$ in \eqref{D_l},
as shown in \eqref{eq:z_v}.
Therefore, recalling that the diffusion algorithm \eqref{diffusion} can be viewed as \eqref{Exact_diffusion} with $\Bx^k-\Bv^k\approx 0$,
we can obtain the following diffusion-like algorithm (CPGD) from \eqref{alg:diffusion_pre} by the similar approximation $\BD\Bx^k - \Bz \approx 0$ for the second line of \eqref{alg:diffusion_pre}:
\begin{align}\label{alg:CPGD_derivation}
    \begin{array}{lll}
 &\Bx^{k+1} = T(\Bx^k - \alpha \nabla_\Bx \hat{f}(\Bx^k)) 
\end{array}
\end{align}
with $T:\BR^d\to\BR^d$ defined as $T(\Bx)= (\BD^\top\BD)^{-1}\BD^\top P_{\Pi_{l\in\QGactive}\CD_l}^\BQ(\BD\Bx)$.
Note that the operator $T$, which will be defined as the \textit{clique-based projection} in Appendix \ref{A:CPGD}, is equal to the doubly stochastic matrix $\BPhi$ in Proposition \ref{prop:double_stochastic} for $\CD_l$ in \eqref{D_l}.

\section{Clique-based projected gradient descent (CPGD) algorithm}\label{A:CPGD}

We here formalize the generalization of the diffusion algorithm (CPGD) in \eqref{alg:CPGD_derivation}.
We provide detailed convergence analysis, which guarantees the exact convergence under diminishing step sizes and an inexact convergence rate over fixed ones.
Moreover, we provide Nesterov's acceleration and an improved convergence rate.

This section highlights the well-behavedness of clique-wise coupling that enables similar theoretical and algorithmic properties to consensus optimization (diffusion algorithm).

\subsection{Clique-based Projected Gradient Descent (CPGD)}

Consider Problem \eqref{problem_clique-wise_constraints_noprox} with closed convex sets $\CD_l\subset \BR^{d^l},\,l\in\QGactive$.
We suppose Assumptions \ref{assumption:sumCl=N}--\ref{assu:prob}.

To this problem, the CPGD is given as follows:
\begin{align}\label{CPGD}
    \Bx^{k+1} = T^p( \Bx^k - \lambda^k \nabla_\Bx \hat{f}(\Bx^k) ),
\end{align}
where $T:\BR^d\to\BR^d$ is the \textit{clique-based projection} for 
\begin{equation}\label{eq:CPGD_D}
    \CD = \bigcap_{l\in\QGactive}\{\Bx\in\BR^d : x_\CCl \in \CD_l \},
\end{equation}
$T^p = \underbrace{T \circ T\circ \cdots \circ T}_{p}$,
$\hat{f}(\Bx) = \sum_{i=1}^n \hat{f}_i(x_i)$, and $\lambda^k$ is a step size.
The clique-based projection $T$ is defined as follows.
\begin{defi}\label{def:clique-based_projection}
Suppose Assumption \ref{assumption:sumCl=N}.
For a non-empty closed convex set $\CD$ in \eqref{eq:CPGD_D}, a graph $\CG$, and its cliques $\CC_l,\,l\in\QGactive$, the \textit{clique-based projection} $T:\BR^d\to\BR^d$ of $x\in\BR^d$ onto $\CD$ is defined as $T(\Bx) = [T_1(x_{\CN_1})^\top,\ldots,T_n(x_{\CN_n})^\top]^\top$
with \begin{equation}
    \label{eq:map_T_i}
    T_i(x_{\CN_i}) = \frac{1}{|\QiGactive|} \sum_{l\in\QiGactive} 
     E_{l,i} P_{\CD_l}^{Q_l} (x_{\CC_l}) 
\end{equation}
for each $i\in \CN$.
\end{defi}

The clique-based projection can be represented as 
    $T(\Bx) = (\BD^\top\BD)^{-1}\BD^\top P^\BQ_{\Pi_{l\in\QGactive}\CD_l}(\BD\Bx).$

The clique-based projection $T$ has many favorable operator-theoretic properties as follows.
\begin{prop}\label{prop:clique-based_projection}
Suppose Assumption \ref{assumption:sumCl=N}.
For the closed convex set $\CD$ in \eqref{eq:CPGD_D} and clique-based projection $T$ in Definition \ref{def:clique-based_projection} onto $\CD$, the following statements hold:
    \begin{itemize}
    \setlength{\itemsep}{0.05cm}
        \item[(a)] The operator $T$ is firmly nonexpansive, i.e., 
        $\|T(\Bx)-T(\Bw)\|^2 \leq  (\Bx-\Bw)^\top (T(\Bx)-T(\Bw))$ 
        holds for any $\Bx,\Bw\in\BR^d$.
        \item[(b)] The fixed points set of $T$ satisfies $\mathrm{Fix}(T)=\CD$.
        \item[(c)] For any $\Bx\in\BR^d\setminus\CD$ and any $\Bw\in\CD$, $\|T(\Bx)-\Bw\|<\|\Bx-\Bw\|$ holds.
        \item[(d)] For any $\Bx\in\BR^d$, $T^\infty(\Bx)=\lim_{p\to \infty} T^p(\Bx) \in \CD$ holds.
    \end{itemize}
\end{prop}
\begin{proof}
See Appendix \ref{proof_prop_cp}.
\end{proof}

The convergence properties of the CPGD over various step sizes are presented as follows.
Note that the CPGD with fixed step sizes does not exactly converge to an optimal solution like the DGD and diffusion methods for consensus optimization.
\begin{theorem}\label{theorem:CPGD_non_acc}
Consider Problem \eqref{problem_clique-wise_constraints} with closed convex sets $\CD_l,\, l\in \QGactive$. 
Consider the CPGD algorithm in \eqref{CPGD}.
Suppose Assumptions \ref{assumption:sumCl=N}--\ref{assu:prob}.
\begin{itemize}
    \item[(a)] Let a positive sequence $\{\lambda^k\}$ satisfy $\lim_{k\to \infty} \lambda^k=0$, $\sum_{k=1}^\infty \lambda^k =\infty$, and $\sum_{k=1}^\infty (\lambda^k)^2 < \infty$.\footnote{For example, $\lambda^k=1/k$ satisfies the conditions.}
    Assume that $\CD$ is bounded.
    Then, for any $\Bx^0\in\BR^d$ and any $p\in\mathbb{N}$, $\Bx^k$ converges to an optimal solution $\Bx^*\in\argmin_{\Bx\in\CD}\hat{f}(\Bx)$.
     \item[(b)] Let a positive sequence $\{\lambda^k\}$ satisfy $\lim_{k\to \infty} \lambda^k=0$, $\sum_{k=1}^\infty \lambda^k =\infty$, and $\sum_{k=1}^\infty |\lambda^k - \lambda^{k+1}| < \infty$.
    \footnote{For example, $\lambda^k=1/k$ and $\lambda^k=1/\sqrt{k}$ satisfy the conditions.}
    Additionally, assume that $\hat{f}(\Bx)$ is strongly convex.
 Then $\Bx^k$ converges to the unique optimal solution $\Bx^*=\argmin_{\Bx\in\CD} \hat{f}(\Bx)$ for any $\Bx^0\in\BR^d$ and any $p\in\mathbb{N}$.
    \item[(c)] Let $\lambda^k=\alpha \in(0,1/\hat{L}]$ for any $k\in\mathbb{N}$.
    Let $J:\BR^d\to\mathbb{R}$ be 
    \begin{equation}\label{eq:J}
    J(\Bx)=\hat{f}(\Bx)+ V(\Bx)/\alpha
    \end{equation}
    with 
    \begin{equation}\label{eq:V}
    V(\Bx) = \frac{1}{2} \sum_{l\in\QGactive} \|x_{\CC_l}-P_{\CD_l}^{Q_l}(x_{\CC_l})\|^2_{Q_l}.
\end{equation}
    Then, for any $\Bx^0\in\BR^d$ and $p=1$,
    \begin{equation}\label{eq:fixed_bound}
         J(\Bx^k)-J(\Bx^*)\leq \frac{\|\Bx^0-\Bx^*\|^2}{2\alpha k} 
    \end{equation}
    holds for $\Bx^*\in \argmin_{\Bx\in\CD} \hat{f}(\Bx)$.
\end{itemize}
\end{theorem}
\begin{proof}(a) From Proposition \ref{prop:clique-based_projection}a-b, the CPGD in \eqref{CPGD} can be regarded as the hybrid steepest descent in \cite{yamada2001hybrid,yamada2002numerically} for any $p\in\mathbb{N}$.
Hence, Theorem \ref{theorem:CPGD_non_acc}a follows from Theorem 2.18, Remark 2.17 in \cite{yamada2002numerically}, and Proposition \ref{prop:clique-based_projection}c. (b) The statement follows from Theorem 2.15 in \cite{yamada2002numerically} and Proposition \ref{prop:clique-based_projection}a-b.
(c) See Appendix \ref{proof_cpgd_acpgd}.
\end{proof}
\begin{rem}
    Using $V$ in \eqref{eq:V}, another expression of the clique-based projection $T$ is obtained as follows.
\begin{prop}\label{prop:V_T}
    Consider the function $V:\BR^d\to\mathbb{R}$ in \eqref{eq:V}.
    Then, it holds for any $\Bx\in\BR^d$ that
\begin{equation}\label{eq:V_T}
    T(\Bx) = \Bx - \nabla_\Bx V(\Bx).
\end{equation}
\end{prop}
 \begin{proof}
     Since each $\CD_l$ is closed and convex, $1/2\,\|x_{\CC_l}-P_{\CD_l}^{Q_l}(x_{\CC_l}) \|_{Q_l}^2$ is differentiable, and thus $V(\Bx)$ in \eqref{eq:V} is also differentiable.
     Then, for all $i\in\CN$, we have
        $\nabla_{x_i} V(\Bx) 
          = \sum_{l\in \QiGactive} \frac{1}{|\QiGactive|} (x_i - E_{l,i} P_{\CD_l}^{Q_l}(x_{\CC_l}))
         = x_i - \frac{1}{|\QiGactive|} \sum_{l\in \QiGactive}  E_{l,i} P_{\CD_l}(x_{\CC_l})
         = x_i -T_i(x_{\CN_i})$
     from \eqref{eq:neighbor_clique} and \eqref{eq:map_T_i}.
     Hence, we obtain \eqref{eq:V_T}.
\end{proof}
  From Proposition \ref{prop:V_T}, we can interpret the CPGD as a variant of the proximal gradient descent \cite{ryu2022large,condat2023proximal,beck2009fast} since the clique-based projection $T$ can be represented as
    $ T(\Bx) = \argmin_{\Bx' \in\BR^d} \: \frac{1}{2}\|\Bx-\Bx'\|^2 + V(\Bx) + \nabla_\Bx V(\Bx)^\top(\Bx'-\Bx).$
    In this sense, the CPGD is a generalization of the conventional projected gradient descent (PGD).
    When $\CG$ is complete, the CPGD equals PGD because $\QG=\{1\}$ and $\CC_1=\CN$ hold for complete graphs.
\end{rem}
\begin{rem}
    A benefit of the CPGD over the CD-DYS is its simple structure which makes its analysis and extension easy.
     We can easily evaluate stochastic and online variants of the CPGD using the same strategy as the online projected gradient descent \cite{hazan2007logarithmic} from Proposition \ref{prop:clique-based_projection}.
\end{rem}

\subsection{Nesterov's acceleration}
    The CPGD with fixed step sizes can be accelerated up to the inexact convergence rate of $O(1/k^2)$ with Nesterov's acceleration \cite{nesterov1983method,beck2009fast}.
    The accelerated CPGD (ACPGD) is given as follows:
\begin{align}\label{ACPGD}
    &\Bx^{k+1} = T^p( \hat{\Bx}^k - \lambda^k \nabla_\Bx \hat{f}(\hat{\Bx}^k) )\nonumber \\
    &\hat{\Bx}^{k+1} = \Bx^{k+1} - \frac{\sigma^k-1}{\sigma^{k+1}} (\Bx^{k+1} - \Bx^{k}),
\end{align}
where $\hat{\Bx}^0 = \Bx^0$ and $\sigma^{k+1}= (1+\sqrt{1+4\sigma^2})/2$ with $\sigma^0=1$.
This algorithm can also be implemented in a distributed manner.

The convergence rate is proved as follows.
\begin{theorem}\label{theorem:CPGD_acc}
    Consider Problem \eqref{problem_clique-wise_constraints} with closed convex sets $\CD_l,\, l\in \QGactive$ and the ACPGD algorithm \eqref{ACPGD}.
    Suppose Assumption \ref{assumption:sumCl=N}.
    Assume that $\CD\subset\BR^d$ in \eqref{eq:CPGD_D} is a non-empty closed convex set.
    Let $p=1$ and $\lambda^k=\alpha \in (0,1/ \hat{L}]$ for all $k$.
    Then, for any initial state $\Bx^0=\hat{\Bx}^0\in\BR^d$, the following inequality holds:
    \begin{equation}\label{eq:nesterov_bound}
        J(\Bx^k) - J(\Bx^*) \leq \frac{2\|\Bx^0-\Bx^*\|^2}{\alpha k^2},
    \end{equation}
     where $\Bx^* \in \argmin_{\Bx\in\CD} \hat{f}(\Bx)$ and $J(\Bx)$ is given as \eqref{eq:J}.
\end{theorem}
\begin{proof}
See Appendix \ref{proof_cpgd_acpgd}.
\end{proof}

\subsection{Proof of Proposition \ref{prop:clique-based_projection}}\label{proof_prop_cp}

As a preliminary, we present important properties of the function $V(\Bx)$ in \eqref{eq:V} for $\CD$ in \eqref{eq:CPGD_D} as follows. Note that the function $V$ in \eqref{eq:V} is convex because of the convexity of each $\CD_l$.
\begin{prop}\label{prop:V}
For $V(\Bx)$ in \eqref{eq:V} and a non-empty closed convex set $\CD$ in \eqref{eq:CPGD_D}, $V(\Bx)=0 \Leftrightarrow \Bx\in \CD$ holds.
\end{prop}
\begin{proof}
    If $V(\Bx)=0$ for $\Bx\in\BR^d$, we obtain $x_{\CC_l}=P_{\CD_l}^{Q_l}(x_{\CC_l}) \in \CD_l$ for all $l\in\QGactive$, which yields $\Bx\in\CD$ because of \eqref{eq:CPGD_D}.
    Conversely, if $\Bx\in\CD$, then we have $x_{\CC_l}\in \CD_l$ for all $l\in\QGactive$.
    Thus, $V(\Bx)=0$ holds.
\end{proof}
\begin{prop}\label{prop:V_smooth}
    The function $V(\Bx)$ in \eqref{eq:V} is a $1$-smooth function, i.e., its gradient $\nabla_\Bx V(\Bx)$ is $1$-Lipschitzian.
\end{prop}
\begin{proof}
    From Definition \ref{def:clique-based_projection}, 
    1-cocoercivity of $P_{\CD_l}^{Q_l}$ (see \cite{bauschke2017correction}),
    and Proposition \ref{prop:V_T}, we obtain the following for any $\Bx,\Bw\in\BR^d$:
    \begin{align*}
        &\|\nabla_\Bx V(\Bx)- \nabla_\Bx V(\Bw)\|^2 = \|(\Bx-\Bw) - (T(\Bx)-T(\Bw))\|^2 \\
        =& \|\Bx-\Bw\|^2 + \|T(\Bx)-T(\Bw)\|^2 -2 (\Bx-\Bw)^\top(T(\Bx)-T(\Bw)) \\
        =&\|\Bx-\Bw\|^2 + \|T(\Bx)-T(\Bw)\|^2 \\
         & -2 \sum_{l\in\QGactive} (x_{\CC_l}-w_{\CC_l})^\top Q_l (P_{\CD_l}^{Q_l}(x_{\CC_l})-P_{\CD_l}^{Q_l}(w_{\CC_l})) \\
        \leq & \|\Bx-\Bw\|^2+\|T(\Bx)-T(\Bw)\|^2\\
         &\quad -2 \sum_{l\in\QGactive}  \|P_{\CD_l}^{Q_l}(x_{\CC_l})-P_{\CD_l}^{Q_l}(w_{\CC_l})\|^2_{Q_l} \\
         \leq & \|\Bx-\Bw\|^2-\|T(\Bx)-T(\Bw)\|^2\leq\|\Bx-\Bw\|^2.
     \end{align*}
     The last line follows from \eqref{eq:T_nonexpansive_proof} in the proof of Proposition \ref{prop:clique-based_projection}a.
     It completes the proof.
\end{proof}

With this in mind, we prove Proposition \ref{prop:clique-based_projection} as follows.

(a) 
From Jensen's inequality and
the quasinonexpansiveness of convex projection operators \cite{bauschke2017correction}, the following inequality holds for any $\Bx,\,\Bw\in\BR^d$:
\begin{align}\label{eq:T_nonexpansive_proof}
    &(T(\Bx)-T(\Bw))^\top(\Bx-\Bw) \nonumber \\ 
    &= \sum_{l\in\QGactive} (x_\CCl-w_\CCl)^\top Q_l (P_{\CD_l}^{Q_l}(x_\CCl)-P_{\CD_l}^{Q_l}(w_\CCl)) \nonumber\\
    &\geq  \sum_{l\in\QGactive} \|P_{\CD_l}^{Q_l}(x_\CCl)-P_{\CD_l}^{Q_l}(w_\CCl)\|_{Q_l}^2\nonumber\\
    &=\sum_{i=1}^n  \frac{1}{|\QiGactive|} \|E_{l,i}P_{\CD_l}^{Q_l}(x_\CCl) - E_{l,i}P_{\CD_l}^{Q_l}(w_\CCl) \|^2\nonumber\\
    &\geq\sum_{i=1}^n \|T_i(x_{\CN_i}) - T_i(w_{\CN_i})\|^2 = \|T(\Bx)-T(\Bw)\|^2.
\end{align}
Thus, we obtain $\|T(\Bx)-T(\Bw)\|^2\leq (T(\Bx)-T(\Bw))^\top(\Bx-\Bw)$. 

(b) $\CD\subset \mathrm{Fix}(T)$ holds because $x_{\CC_l}=P_{\CD_l}^{Q_l}(x_{\CC_l})$ holds for any $\Bx\in\CD$ and all $l\in\QGactive$.
In the following, we prove the converse inclusion $\mathrm{Fix}(T)\subset\CD$.
Let $\Bw\in\CD$. Then, it suffices to show $\hat{\Bw}\in \mathrm{Fix}(T)\setminus\{\Bw\}  \Rightarrow \hat{\Bw}\in \CD$.
From $\hat{\Bw}\in\mathrm{Fix}(T)$, we obtain $\hat{w}_i=T_i(\hat{w}_{\CN_i})$ for all $i\in\CN$.
In addition, from Jensen's inequality and the quasinonexpansiveness of convex projection operators \cite{bauschke2017correction}, we have
\begin{align*} 
    &\|\Bw-\hat{\Bw}\|^2
    \geq \sum_{l\in\QGactive} \|w_{\CC_l} - P_{\CD_{\CC_l}}^{Q_l}(\hat{w}_\CC)\|^2_{Q_l} \nonumber\\
    &=\sum_{i=1}^n \sum_{l\in\QiGactive} \frac{1}{|\QiGactive|} \|w_i-E_{l,i} P_{\CD_l} (\hat{w}_{\CC_l})\|^2  \nonumber \\
   &\geq \sum_{i=1}^n  \| w_i - 
   \underbrace{\sum_{l\in \QiGactive} \frac{1}{|\QiGactive|} E_{l,i} P_{\CD_l} (\hat{w}_{\CC_l})}_{=T_i(\hat{w}_{\CN_i})=\hat{w}_{\CN_i}} \|^2
   = \|\Bw-\hat{\Bw}\|^2.
\end{align*}
Thus, from the equality condition of Jensen's inequality, we obtain $w_i-E_{l,i} P_{\CD_k} (\hat{w}_{\CC_k}) =w_i-E_{l,i} P_{\CD_l} (\hat{w}_{\CC_l}) $ for all $\CC_k,\CC_l\,(k,l\in\QiGactive)$ for all $i\in\CN$.
Then, we have $E_{l,i} P_{\CD_k} (\hat{w}_{\CC_k}) =E_{l,i} P_{\CD_l} (\hat{w}_{\CC_l}) $ for all $\CC_k,\CC_l\,(k,l\in\QiGactive)$.
Therefore, since $\hat{\Bw}\in\mathrm{Fix}(T)$, we have $2V(\hat{\Bw})=\sum_{i=1}^n \sum_{l\in \QiGactive} \frac{1}{|\QiGactive|} \|\hat{w}_i-E_{l,i} P_{\CD_l} (\hat{w}_{\CC_l})\|^2= \sum_{i=1}^n \| \hat{w}_i - T_i(\hat{w}_{\CN_i})\|^2=0$.
Thus, $\hat{\Bw}\in\CD$ holds from Proposition \ref{prop:V}.

(c) For a non-empty closed convex set $\CD$ in \eqref{eq:CPGD_D} and $\Bx\in\BR^d\setminus\CD$, there exists $\hl\in\QGactive$ such that $\|x_{\CC_{\hl}}-P_{\CD_{\hl}}(x_{\CC_{\hl}})\|_{Q_{\hl}}>0$.
Hence, for $\hl\in\QGactive$, $\Bx\in\BR^d\setminus\CD$, and $\Bw \in\CD$, we have $\|x_{\CC_{\hl}}-w_{\CC_{\hl}}\|^2_{Q_{\hl}}>\|P_{\CD_{\hl}}^{Q_{\hl}}(x_{\CC_{\hl}})-w_{\CC_{\hl}}\|^2_{Q_{\hl}}$ because 
\begin{align*}
    &\|x_{\CC_l}-z_{\CC_l}\|^2_{\mathrm{diag}(\gamma_{\CC_l})} \\
    =& \|x_{\CC_l}-P_{\CD_l}(x_{\CC_l})\|^2_{\mathrm{diag}(\gamma_{\CC_l})} + \|P_{\CD_l}(x_{\CC_l})-z_{\CC_l}\|^2_{\mathrm{diag}(\gamma_{\CC_l})} \\
    & \quad -2 (x_{\CC_l}-P_{\CD_l}(x_{\CC_l}))^\top \mathrm{diag}(\gamma_{\CC_l})(z_{\CC_l}-P_{\CD_l}(x_{\CC_l})) \\
    >& \|P_{\CD_l}(x_{\CC_l})-z_{\CC_l}\|^2_{\mathrm{diag}(\gamma_{\CC_l})},
\end{align*}
where the last line follows from the projection theorem (see Theorem 3.16 in \cite{bauschke2017correction}).
Thus, by Jensen's inequality and the nonexpansiveness of $P_{\CD_l}^{Q_l}$ \cite{bauschke2017correction}, for any $\Bx\in\BR^d\setminus\CD$ and $\Bw\in\CD$,
we obtain $ \|\Bx-\Bw\|^2 = \sum_{l\in \QGactive} \|x_{\CC_l}-w_{\CC_l}\|^2_{Q_l}  
     > \sum_{l\in \QGactive} \|P_{\CD_l}^{Q_l}(x_{\CC_l})-w_{\CC_l}\|^2_{Q_l}  
    \geq \sum_{i=1}^n \|
    \sum_{l\in\QiGactive} 
    \frac{1}{|\QiGactive|}E_{l,i}P_{\CD_l}^{Q_l} (x_{\CC_l})\|^2
    = \|T(\Bx)-\Bw\|^2.$
Hence, $\|T(\Bx)-\Bw\|<\|\Bx-\Bw\|$ for any $\Bx\in\BR^d\setminus\CD$ and $\Bw\in\CD$.

(d) For $\Bx\in\BR^d$, we define $\{a_k\}$ as $a_{k+1} = T(a_k)$ with $a_0=\Bx$.
Then, we obtain $\lim_{k\to\infty} a_{k+1} = \lim_{k\to\infty} T(a_k)$.
Thus, from the continuity of $T$ shown in Proposition \ref{prop:clique-based_projection}a, we have $T^\infty(x)= \lim_{k\to\infty} a_{k+1} =T( \lim_{k\to\infty} a_k) = T(T^\infty(x))$.
Hence, Proposition \ref{prop:clique-based_projection}b yields $T^\infty (x)\in\mathrm{Fix}(T) = \CD$.

\subsection{Proof of Theorems \ref{theorem:CPGD_non_acc}c and \ref{theorem:CPGD_acc}}\label{proof_cpgd_acpgd}

Here, we show the proofs of Theorems \ref{theorem:CPGD_non_acc}c and \ref{theorem:CPGD_acc}.
These proofs are based on the convergence theorems for the ISTA and FISTA (Theorems 3.1 and 4.4 in \cite{beck2009fast}), respectively.
\paragraph{Supporting Lemmas}
Before proceeding to prove the theorems, we show some inequalities corresponding to those obtained from Lemma 2.3 in \cite{beck2009fast}, which is a key to proving the convergence theorems.
Note that a differentiable function $h:\mathbb{R}^m\to\mathbb{R}$ is convex if and only if
\begin{equation}\label{eq:convex_inquality}
    h(\Bw)\geq h(\Bx) + \nabla h(\Bx)^\top(\Bw-\Bx)
\end{equation}
holds for any $\Bx,\Bw\in\BR^d$.
If $h$ is $\beta$-smooth and convex,
\begin{align}
\label{eq:smooth1}
    h(\Bw) &\leq h(\Bx) + \nabla h(\Bx)^\top(\Bw-\Bx) + \frac{\beta}{2} \|\Bw-\Bx\|^2 \\
\label{eq:smooth2}
    h(\Bw) &\geq h(\Bx) + \nabla h(\Bx)^\top(\Bw-\Bx) + \frac{1}{2\beta} \|\nabla h(\Bx)-\nabla h(\Bw)\|^2
\end{align}
hold for any $\Bx,\Bw\in\BR^d$. For details, see textbooks on convex theory, e.g.,
Theorem 18.15 in \cite{bauschke2017correction}.

In preparation for showing lemmas, let $\alpha \in (0,1/\hat{L}]$ and 
$V_\alpha (\Bx) = V(\Bx)/\alpha$
with $V(\Bx)$ in \eqref{eq:V}.
Additionally, for $\Bs\in\BR^d$, we define $\hat{F}_\Bw :\BR^d\to \mathbb{R}$ with some $\Bw\in\BR^d$ as
\begin{equation}\label{def:F_y}
 \hat{F}_\Bw (\Bs) = \hat{f}(\Bs) + V_\alpha(\Bw) + \nabla_\Bx V_\alpha(\Bw)^\top (\Bs-\Bw).
\end{equation}
For $\hat{F}_\Bw (\Bs)$ in \eqref{def:F_y}, the following inequalities hold.
\begin{prop}
Assume that $\hat{f}$ is $\hat{L}$-smooth and convex.
    Let $\Bw = \Bx-\alpha\nabla_\Bx \hat{f}(\Bx)$.
    Then, 
    \begin{equation}\label{eq:lem_acc_1}
        \hat{F}_\Bw(T(\Bw)) \leq \hat{F}_\Bw (\Bxi) + \frac{1}{\alpha} (\Bx-T(\Bw))^\top (\Bx-\Bxi) - \frac{1}{2\alpha} \|\Bx-T(\Bw)\|^2
    \end{equation}
    holds for any $\Bxi\in\BR^d$.
\end{prop}
\begin{proof}
    Let $G_\Bw(\Bs) = \hat{f}(\Bs) + \nabla_\Bx V_\alpha(\Bw)^\top (\Bs-\Bw) $ and $\Bxi\in\BR^d$. Then, by using $\hat{L}$-smoothness of $\hat{f}$, $\nabla_\Bx \hat{f}(\Bx) = (\Bx-\Bw)/\alpha$, and $\nabla_\Bx V_\alpha (\Bw) = (\Bw-T(\Bw))/\alpha$ (see Proposition \ref{prop:V_T}), 
    \fontsize{9.8pt}{2.5pt}\selectfont
    \begin{align*}
        &G_\Bw(T(\Bw)) 
        = \hat{f}(T(\Bw)) + \nabla_\Bx V_\alpha(\Bw)^\top (T(\Bw)-\Bw) \\
        \leq& \hat{f}(\Bx) - \nabla_\Bx \hat{f}(\Bx)^\top (\Bx-T(\Bw)) + \frac{1}{2\alpha} \|\Bx-T(\Bw)\|\\
        &+ \nabla_\Bx V_\alpha(\Bw)^\top (T(\Bw)-\Bw) \\
        \leq& \hat{f}(\Bxi) + \nabla_\Bx \hat{f}(\Bx)^\top (\Bx-\Bxi) - \nabla_\Bx \hat{f}(\Bx)^\top(\Bx-T(\Bw)) \\
        +& \frac{1}{2\alpha} \|\Bx-T(\Bw)\|^2 
        + \nabla_\Bx V_\alpha(\Bw)^\top \underbrace{(T(\Bw)-\Bw)}_{=(\Bxi-\Bw)+(T(\Bw)-\Bxi)} 
         \\
        =& G_\Bw(\Bxi) + \frac{1}{\alpha}(\Bx-T(\Bw))(T(\Bw)-\Bxi) + \frac{1}{2\alpha}\|\Bx-T(\Bw)\|^2\\
        =&  G_\Bw(\Bxi) + \frac{1}{\alpha}(\Bx-T(\Bw))^\top (\Bx-\Bxi) -\frac{1}{2\alpha} \|\Bx-T(\Bw)\|^2
        \end{align*}
        \normalsize
        is obtained from \eqref{eq:convex_inquality} and \eqref{eq:smooth1}.
        Thus, adding $V_\alpha(\Bw)$ to both sides, we obtain \eqref{eq:lem_acc_1}. 
\end{proof}

\begin{prop}
    Let $\Bx^{k+1} = T(\Bw^k)$ with some $\{\Bw^k\}\subset \BR^d$. Then, it holds that
    \begin{align}\label{eq:lem_acc_2}
        \hat{F}_{\Bw^k} (\Bx^k) + \frac{\alpha}{2} \|\nabla_\Bx V_\alpha(\Bw^k))\|^2 
        \leq
        \hat{F}_{\Bw^{k-1}} (\Bx^k) + \frac{\alpha}{2} \|\nabla_\Bx V_\alpha(\Bw^{k-1})\|^2.
    \end{align}
\end{prop}
\begin{proof}
By $1/\alpha$-smoothness of $V_\alpha(\Bx)$ (see Proposition \ref{prop:V_smooth}) and Proposition \ref{prop:V_T},
\begin{align*}
    &\hat{F}_{\Bw^{k-1}} (\Bx^k) = \hat{f}(\Bx^k) + V_\alpha(\Bw^{k-1}) 
    + \nabla_\Bx V_\alpha(\Bw^{k-1})^\top (\Bx^k-\Bw^{k-1}) \\
     =& \hat{f}(\Bx^k) + V_\alpha(\Bw^{k-1}) - \alpha\|\nabla_\Bx V_\alpha(\Bw^{k-1})\|^2 \\
    \geq& \hat{f}(\Bx^k) + V_\alpha(\Bw^k) + \nabla_\Bx V_\alpha(\Bw^k)^\top (\Bw^{k-1}-\Bw^k) \\
    &+ \frac{\alpha}{2} \|\nabla_\Bx V_\alpha(\Bw^{k-1})-\nabla_\Bx V_\alpha(\Bw^k)\|^2 - \alpha\|\nabla_\Bx V_\alpha(\Bw^{k-1})\|^2 \\
    =& \hat{f}(\Bx^k) + V_\alpha(\Bw^k) + \nabla_\Bx V_\alpha(\Bw^{k})^\top (\Bx^k-\Bw^k) \\
    &+ \nabla_\Bx V_\alpha(\Bw^k)^\top (\Bw^{k-1}-\Bx^k)
    \\ 
    &+ \frac{\alpha}{2} \|\nabla_\Bx V_\alpha(\Bw^{k-1})-\nabla_\Bx V_\alpha(\Bw^k)\|^2 - \alpha\|\nabla_\Bx V_\alpha(\Bw^{k-1})\|^2 \\
    =& \hat{F}_{\Bw^k}(\Bx^k) + \frac{\alpha}{2} \|\nabla_\Bx V_\alpha(\Bw^k)\|^2 -\frac{\alpha}{2} \|\nabla_\Bx V_\alpha(\Bw^{k-1})\|^2
\end{align*}
    is obtained from \eqref{eq:smooth2}. Hence, \eqref{eq:lem_acc_2} holds. 
\end{proof}

With this in mind, we consider the following update rule with $\hat{\Bx}(0) = \Bx(0)$ and some 
$\{\theta^k\}\subset \mathbb{R}$:
\begin{align}\label{nesterov_theta}
    \Bw^k &= \hat{\Bx}^k - \alpha \nabla_\Bx \hat{f} (\hat{\Bx}(k)) \nonumber \\
        \Bx^{k+1}  &= T(\Bw^k)  \nonumber \\
    \hat{\Bx}^{k+1} &=  \Bx^{k+1} + \theta^k (\Bx^{k+1}-\Bx^k).
\end{align}
In addition, we define $\Theta^k:\BR^d\to\mathbb{R}$ as 
\begin{equation}\label{def:H}
    \Theta^k = \hat{F}_{\Bw^{k-1}}(\Bx^k) + \frac{\alpha}{2}\|V_\alpha(\Bw^{k-1})\|^2
\end{equation}
with $\hat{F}_\Bw$ in \eqref{def:F_y}.
By $\Bx^k-\Bw^{k-1} = -\alpha \nabla_\Bx V_\alpha(\Bw^{k-1})$, 
$\Theta^k$ can be rewritten as $\Theta^k = \hat{f}(\Bx^k)+V_\alpha(\Bw^{k-1})- \frac{1}{2\alpha} \|\Bw^{k-1}-T(\Bw^{k-1})\|^2 = \hat{f}(\Bx^k)+V_\alpha(\Bw^{k-1}) - \frac{1}{2\alpha} \|\Bw^{k-1}-\Bx^{k}\|^2$.

Remarkably, $\Theta^k$ in \eqref{def:H} satisfies the following lemma.
\begin{lemma}
    Consider the sequence generated by \eqref{nesterov_theta}. Then,
    \begin{equation}\label{eq:bound_H_k}
    J(\Bx^k)=\hat{f}(\Bx^k) + V_\alpha(\Bx^k) \leq \Theta^k.
    \end{equation}
\end{lemma}
\begin{proof}
    In light of $1/\alpha$-smoothness of $V_\alpha$ and $\nabla_\Bx V_\alpha(\Bw^{k-1})= -(\Bw^{k-1}-\Bx^k)/\alpha$, we obtain $V_\alpha(\Bx^k) 
        \leq V_\alpha(\Bw^{k-1}) 
        + \nabla_\Bx V_\alpha(\Bw^{k-1})^\top (\Bw^{k-1}-\Bx^{k})  + \frac{1}{2\alpha} \|\Bw^k-\Bx^{k}\|^2 
        = V_\alpha(\Bw^{k-1}) - \frac{1}{2\alpha} \|\Bw^{k}-\Bx^{k}\|^2.$
    Hence, adding $\hat{f}(\Bx^k)$ to both sides yields \eqref{eq:bound_H_k}. 
\end{proof}
Furthermore,
the following inequality holds. This is essential to Theorem \ref{theorem:CPGD_non_acc}c and \ref{theorem:CPGD_acc}.
\begin{lemma}
    For the sequence generated by \eqref{nesterov_theta} and $\Theta^k$ defined in \eqref{def:H}, it holds that
    \begin{align}\label{eq:lem_H_k}
    \small
        \Theta^k-\Theta^{k+1} &\geq \frac{1}{2\alpha} \| \hat{\Bx}^k- \Bx^{k+1} \|^2  
        + \frac{1}{\alpha} (\Bx^{k+1}-\hat{\Bx}^k)^\top (\hat{\Bx}^k-\Bx^k).
    \end{align}
\end{lemma}
\begin{proof}
    Substituting $\Bx=\Bx^{k+1},\,\Bw=\Bw^k,$ and $\Bxi=\Bx^k$ into \eqref{eq:lem_acc_1}, we obtain
    \fontsize{9.8pt}{2.5pt}\selectfont
    \begin{align*}
        &\Theta^{k+1} = \hat{f}(\Bx^{k+1}) + V_\alpha(\Bw^k) \\
        &+ \nabla_\Bx V_\alpha(\Bw^k)^\top (\Bx^{k+1}-\Bw^{k}) + \frac{\alpha}{2}\|\nabla_\Bx V_\alpha(\Bw^k)\|^2 \\
        &\leq \hat{f}(\Bx^k) + V_\alpha(\Bw^k) \\
        &+ \nabla_\Bx V_\alpha(\Bw^k)^\top (\Bx^k-\Bw^k) + \frac{\alpha}{2} \|\nabla_\Bx V_\alpha(\Bw^k)\|^2 \\
        &+ \frac{1}{\alpha} (\hat{\Bx}^k-\Bx^{k+1})^\top (\hat{\Bx}^k-\Bx^k) - \frac{1}{2\alpha} \|\hat{\Bx}^k-\Bx^{k+1}\|^2 \\
        &= F_{\Bw^k}(\Bx^k) + \frac{\alpha}{2}\|\nabla_\Bx V_\alpha(\Bw^k)\|^2 \\
        &+  \frac{1}{\alpha} (\hat{\Bx}^k-\Bx^{k+1})^\top (\hat{\Bx}^k-\Bx^k) 
        -\frac{1}{2\alpha} \|\hat{\Bx}^k-\Bx^{k+1}\|^2 \\
        &\leq F_{\Bw^{k-1}}(\Bx^k) + \frac{\alpha}{2}\|\nabla_\Bx V_\alpha(\Bw^{k-1})\|^2 \\
        &+  \frac{1}{\alpha} (\hat{\Bx}^k-\Bx^{k+1})^\top (\hat{\Bx}^k-\Bx^k) 
        -\frac{1}{2\alpha} \|\hat{\Bx}^k-\Bx^{k+1}\|^2 \\
        &= \Theta^k +  \frac{1}{\alpha} (\hat{\Bx}^k-\Bx^{k+1})^\top (\hat{\Bx}^k-\Bx^k) 
        -\frac{1}{2\alpha} \|\hat{\Bx}^k-\Bx^{k+1}\|^2
    \end{align*}
    from \eqref{eq:convex_inquality}, \eqref{eq:smooth1}, and \eqref{eq:lem_acc_2}. Thus, \eqref{eq:lem_H_k} holds. 
\end{proof}
\normalsize

For $\Bx^k$ and an optimal $\Bx^*$, we present the following lemma.
\begin{lemma}
    For $\Bx^*\in \argmin_{\Bx\in\CD} \hat{f}(\Bx)$, it holds that
    \begin{align}\label{eq:lem_H_k_2}
    \hat{f}(\Bx^*)+V_\alpha(\Bx^*) &- \Theta^{k+1} \geq
    \frac{1}{2\alpha} \|\hat{\Bx}^k-\Bx^{k+1}\|^2
    + \frac{1}{\alpha} (\Bx^{k+1}-\hat{\Bx}^k)^\top (\hat{\Bx}^k-\Bx^*).
    \end{align}
\end{lemma}
\begin{proof}
Recalling \eqref{nesterov_theta}, $\hat{L}$-smoothness of $\hat{f}$, and $1/\alpha$-smoothness of $V_\alpha$ for $\alpha\in(0,1/\hat{L}]$, we obtain
\fontsize{9.8pt}{2.5pt}\selectfont
\begin{align*}
    &\Theta^{k+1} \leq \hat{f}(\hat{\Bx}^k) - \nabla_\Bx \hat{f}(\hat{\Bx}^k)^\top (\hat{\Bx}^k-\Bx^{k+1})\\
    & + \frac{1}{2\alpha} \|\hat{\Bx}^k-\Bx^{k+1}\|^2 + V_\alpha(\Bw^k) - \frac{1}{2\alpha} \|\Bw^k-T(\Bw^k)\|^2 \\
    &\leq \hat{f}(\Bx^*) + \nabla_\Bx \hat{f}(\hat{\Bx}^k)^\top (\hat{\Bx}-\Bx^*)- \nabla_\Bx \hat{f}(\hat{\Bx}^k)^\top (\hat{\Bx}-T(\Bw^k))  \\
    & + \frac{1}{2\alpha} \|\hat{\Bx}^k-T(\Bw^k)\|^2+V_\alpha(\Bx^*) - \frac{1}{2\alpha} \|\Bw^k-T(\Bw^k)\|^2 \\
    & + \frac{1}{\alpha}(\Bw^k-T(\Bw^k)^\top (T(\Bw^k)-\Bx^*+\Bw^k-T(\Bw^k)) \\
    &-\frac{1}{2\alpha} \|\Bw^k-T(\Bw^k) - (\Bx^*-T(\Bx^*))\|^2  \\
    &=\hat{f}(\Bx^*) + V_\alpha(\Bx^*) + \frac{1}{\alpha} (\hat{\Bx}^k-\Bx^{k+1})^\top (\hat{\Bx}^k-\Bx^*) 
- \frac{1}{2\alpha} \|\hat{\Bx}^k-\Bx^{k+1}\|^2
\end{align*}
\normalsize
from \eqref{eq:convex_inquality}, \eqref{eq:smooth1}, and \eqref{eq:smooth2},
where the last line is obtained because $\Bx^*=T(\Bx^*)$ holds for $\Bx^*\in\CD$.
Therefore, \eqref{eq:lem_H_k_2} is obtained.
\end{proof}

\paragraph{Proof of Theorem \ref{theorem:CPGD_non_acc}c}
In this proof, assume that $\theta^k=0$ for all $k$.
Then, $\hat{\Bx}^k=\Bx^k$ holds and the algorithm in \eqref{nesterov_theta} equals to the CPGD with $\lambda^k=\alpha\in(0,1/\hat{L}]$ for all $k\in\mathbb{N}$.

In light of \eqref{eq:lem_H_k_2} and $\hat{\Bx}^k=\Bx^k$, we obtain $2\alpha(\Theta^{k+1}-(\hat{f}(\Bx^*)+V_\alpha(\Bx^*)))\leq\|\Bx^*-\Bx^k\|^2$ because $2\alpha(\Theta^{k+1}-(\hat{f}(\Bx^*)+V_\alpha(\Bx^*))) \leq  2(\Bx^k-\Bx^{k+1})^\top (\Bx^k-\Bx^*)
- \frac{1}{2\alpha} \|\Bx^k-\Bx^{k+1}\|^2 
= \|\Bx^*-\Bx^k\|^2 - \|\Bx^*-\Bx^{k+1}\|^2 \leq\|\Bx^*-\Bx^k\|^2.$
Besides, invoking \eqref{eq:lem_H_k}, we have
\begin{equation*}
    2\alpha(\Theta^{k+1}-\Theta^k )\leq \|\Bx^k-\Bx^{k+1}\|^2 \leq 0.
\end{equation*}
Then, following the same procedure as Theorem 3.1 in \cite{beck2009fast}
and using \eqref{eq:bound_H_k}, we obtain \eqref{eq:fixed_bound}. 

\paragraph{Proof of Theorem \ref{theorem:CPGD_acc}}\label{subsec:proof_acc}
Substituting $\theta^k=(\sigma^k-1)/\sigma^{k+1}$ into \eqref{nesterov_theta} yields the ACPGD in \eqref{ACPGD}.

Now, by \eqref{eq:lem_H_k}, \eqref{eq:lem_H_k_2}, and $(\sigma^{k-1})^2 = \sigma^k(\sigma^k-1) $, following the procedure of the proof of Theorem 4.4 in \cite{beck2009fast} gives
\begin{align*}
\small
    &(\sigma^{k-1})^2(\Theta^{k}-J(\Bx^*))
    - (\sigma^{k})^2(\Theta^{k+1}-J(\Bx^*)
    )
    \\
    &
    \leq \frac{1}{2\alpha} (\|\Bzeta^{k+1}\|^2-\|\Bzeta^k\|^2)
\end{align*}
with $J$ in \eqref{eq:J} and $\Bzeta^k=\sigma_k (\hat{\Bx}^k-\Bx^*) - (\sigma^k-1)(\Bx^k-\Bx^*)$.
Thus, summing both sides over $k=1,2,\ldots$ yields
\begin{equation*}
  (\sigma^k)^2(  \Theta^{k+1} - J(\Bx^*) )
  \leq \frac{1}{2\alpha}\|\Bzeta^0\|^2 =\frac{1}{2\alpha}\|\Bx^0-\Bx^*\|^2.
\end{equation*}
By $\sigma^k\geq (k+1)/2$, which can be shown by mathematical induction, we obtain
\begin{equation*}
    \Theta^{k+1} - J(\Bx^*)
    \leq \frac{2\|\Bx^0-\Bx^*\|^2}{\alpha(k+1)^2}. 
\end{equation*}
Therefore, the inequality \eqref{eq:nesterov_bound} follows from \eqref{eq:bound_H_k}.

\end{document}